
%
%
%
\documentstyle{amsppt}
%
%
%
%
\def\R{\Bbb R }

\def\Z{\Bbb Z }
\def\C{\Bbb C }
\loadeurm
\define\k{{\eurm k}}
\define\ce{{\eurm c}}
\define\ka{{\eurm k}}
\define\npos{{\eurm p}}

\define\nl{\hfil\newline}
\define\ldot{\,.\,}
\redefine\l{\lambda}

\define\a{\alpha}
\redefine\d{\delta}
\define\w{\omega}
\define\ep{\epsilon}

\redefine\t{\tau}
\redefine\i{{\,\text{{\rm i}}\,}}
\define\ga{\gamma}
\define\PL{Phys\. Lett\. B}

\define\LMP{Lett\. Math\. Phys\. }
\redefine\CMP{Commun\. Math\. Phys\. }
\define\JMP{J\.  Math\. Phys\. }

\define\FA{Funktional Anal\. i\. Prilozhen\.}
\def\Pnas{Proc\. Natl\. Acad\. Sci\. USA}
\def\PAMS{Proc\. Amer\. Math\. Soc\.}
\define\mapleft#1{\smash{\mathop{\longleftarrow}\limits^{#1}}}
\define\mapright#1{\smash{\mathop{\longrightarrow}\limits^{#1}}}
\define\mapdown#1{\Big\downarrow\rlap{
   $\vcenter{\hbox{$\scriptstyle#1$}}$}}
\define\mapup#1{\Big\uparrow\rlap{
   $\vcenter{\hbox{$\scriptstyle#1$}}$}}
\define\cint #1{\frac 1{2\pi\i}\int_{C_{#1}}}
\define\cintt{\frac 1{2\pi\i}\int_{C_{\tau}}}
\define\cinttp{\frac 1{2\pi\i}\int_{C_{\tau'}}}

\define\cintb{\frac 1{2\pi\i}\int_{Z_{\alpha}}}
\define\cinttt{\frac 1{24\pi\i}\int_{C_{\tau}}}
\define\cinttb{\frac 1{24\pi\i}\int_{Z_{\alpha}}}
\define\cintd{\frac 1{(2\pi \i)^2}\iint\limits_{C_{\tau}\,C_{\tau'}}}

\define\pem{\psi(\ep n)}
\define\im{\text{Im\kern1.0pt }}
\define\re{\text{Re\kern1.0pt }}
\define\res{\text{res}}
\redefine\deg{\operatornamewithlimits{deg}}
\define\ord{\operatorname{ord}}
\define\rank{\operatorname{rank}}
\define\fpz{\frac {d }{dz}}
\define\dzl{\,{dz}^\l}
\define\pfz#1{\frac {d#1}{dz}}
\define\KN {Krichever-Novikov}
\define\KM {Kac-Moody}
\define\Fl{\Cal F^\lambda}
\define\Fn#1{\Cal F^{#1}}
\define\A{\Cal A}
\redefine\L{\Cal L}
\define\G{\Cal G}
\define\Ah{\widehat{\Cal A}}
\define\Lh{\widehat{\Cal L}}
\define\Gh{\widehat{\Cal G}}
\define\Ghe{\widehat{\Cal G}^e}
\define\Ghd{\widehat{\Cal G}^d}
\redefine\g{\frak g}
\define\HG{H_1(\Gamma^*,\Z)}
\define\Tld{...}
\define\Ksig{K(\Sigma)}
\define\ea{x}
\define\eps{\varepsilon}    
\define \fg{{\frak g}}     
\define \nord #1{:\mkern-5mu{#1}\mkern-5mu:}
\magnification=1200
\vsize=21.5truecm
\hsize=16truecm
\hoffset=0.5cm\voffset=0.2cm
\baselineskip=15pt plus 0.2pt
\parskip=8pt
\NoBlackBoxes
\TagsOnRight
\hfill Mannheimer Manuskripte 201

\hfill q-alg/9512016
\vskip 1cm
\topmatter
\title
Sugawara Construction and Casimir Operators for
              Krichever-Novikov Algebras
\endtitle
\rightheadtext{Sugawara Construction}
\leftheadtext{M. Schlichenmaier, O.K. Sheinman}
\author Martin Schlichenmaier
\footnote"*"{Partially~supported~by~the~RiP
 program~of~the~Volkswagen-Stiftung.\hphantom{Das muss aber jetzt reichen}
\hfill }
 and Oleg K. Sheinman\footnotemark"*"
\endauthor
\address
Martin Schlichenmaier,
Department of Mathematics and Computer Science,
University of Mannheim
D-68131 Mannheim, Germany
\endaddress
\email
schlichenmaier\@math.uni-mannheim.de
\endemail
\address
Oleg K. Sheinman,
Independent Moscow University,
Ul. Marshala Biryuzova 4, kor.1, kv.65
Moscow, Russia, 123298
\endaddress
\email
sheinman\@landau.ac.ru
\endemail
\date December 95
\enddate
\keywords
Sugawara construction, highest weight representations, infinite-dimensional
Lie algebras, central extensions, current algebras, Casimir operators
\endkeywords
\subjclass
17B66, 17B67, 17B90, 30F30, 14H55, 81R10, 81T40
\endsubjclass
\abstract
We show how to obtain from highest weight representations of
 Krichever-Novikov algebras of affine type (also
called higher genus affine Kac-Moody algebras)
representations of centrally extended Krichever-Novikov vector field
algebras via the Sugawara construction.
This generalizes  classical results where one obtains
representations of the Virasoro algebra.
Relations between the weights of the corresponding representations
are given and Casimir operators are constructed.
In an appendix the Sugawara construction for the multi-point
situation is done.
\endabstract
\endtopmatter
%
%
%
%
\newcount\refCount
\def\newref#1 {\advance\refCount by 1
\expandafter\edef\csname#1\endcsname{\the\refCount}}

\newref BRRW 
\newref Bor 
\newref BREM 
\newref BREMa 
\newref FLM 
\newref JKL 
\newref KIL 
\newref KACRAI %
\newref	KPC 
\newref	KNFA 
\newref	KNJGP 
\newref	KNPMS 
\newref SADOV 
\newref SCHLRS 
\newref SCHLL 
\newref SCHLTH 
\newref SCHLCT 
\newref SHEIN 
\newref SHEINW 
\newref SHEINNS  
%
%
%
%

\define\knset{2}
\define\knsu{3}
\define\knpr{4}
\define\knwei{5}

%
%
%
%
\document
\vskip 1cm
\head
1. Introduction
\endhead
\def\kn{1}
\vskip 0.2cm
The  Sugawara construction is one of the basic constructions in
two-dimensional conformal field theory.
Also from the point of view of representation theory
it  is at least  for two reasons of importance.
First, it
provides a realization of highest weight representations of
Virasoro-type algebras and
second,  the same technique  is used for
the investigation of Casimir operators for affine type algebras.
Nevertheless, in the mathematical literature there is no consecutive
presentation of this construction for higher genera with all necessary
proofs. We hope that this article will fill this gap. Almost all
necessary ideas are contained in Kac \cite\KACRAI\ (see also the
references therein), Krichever and Novikov \cite\KNFA\ and Bonora et al.
\cite\BRRW\ but each of these three treatments is too restrictive with
respect to the genus, or to the properties of the underlying
 finite dimensional
Lie algebra, or at least to the completeness of the presentation. We give
in this article a survey of the basic results in form
 of a detailed proof.
In addition, we demonstrate the application of this technique in the
construction of  Casimir operators for higher genera and
generalize it to the case of Riemann surfaces with many
punctures.

In the framework of conformal field theory representations of affine
Kac-Moody algebras  and of the Virasoro algebra play a fundamental
role.
These algebras and their representations are also from the
mathematical point of view of great interest.
They give interesting examples of infinite dimensional Lie algebras
which one still can handle. For their representations one has
developed a fairly complete but nevertheless nontrivial
structure theory.
Another very fascinating aspect of them is that they have been proven
to be very useful in proving deep mathematical results.
Let us mention here the explanation of the ``monster and moonshine''
which relates the dimensions of the  irreducible representations
of the monster group  with the coefficients of the $q$-expansion
series for the elliptic modular function $j(q)$, observed by Conway
and explained by Frenkel, Lepowski and Meurman \cite\FLM\
(see also Borcherds \cite\Bor).
Recall that the monster is the largest  exceptional finite simple group.
There are lots of other important applications, like in
the theory of integrable systems and so on.

Let us recall the definitions of these algebras.
For details see \cite\KIL, or \cite\KACRAI\  for a more
pedagogical treatment.
Let $\g$ be a finite dimensional Lie algebra over $\C$.
The current algebra (or loop algebra)  is defined
to be $\G=\g\otimes\C[z,z^{-1}]$ as vector space. It is generated by the
elements $\ x\otimes z^n\ $  with $n\in\Z$, $x\in\g$ and
Lie structure given by
$$[x\otimes z^n\,,\,y\otimes z^m]=[x\,,\,y]\otimes z^{n+m}\ .\tag\kn-1$$
If $\g$ admits an  invariant symmetric non-degenerate
bilinear form $\ (..|..)\ $  then the affine (untwisted)
Kac-Moody algebra is the
centrally extended algebra $\Gh$ of the current algebra.
As vector space $\Gh=\G\oplus\C\, t$ with Lie structure
$$[x(n)\,,\,y(m)]=[x\,,\,y](n+m)+(x|y)\cdot  n\cdot \delta^n_{-m}\cdot t,
\qquad [\,t\,,\,\Gh]=0\ ,\tag\kn-2$$
where we used the usual notation $\ x(n):=x\otimes z^n$.
Special cases are  the situation  where $\g$ is simple and one
takes the  Cartan-Killing form or where $\g$ is abelian and
one takes an arbitrary non-degenerate symmetric bilinear form.
In the latter case  one calls $\Gh$ also a Heisenberg algebra.

Often it is convenient to adjoin an additional element $d$, a
derivation,  to
$\Gh$ and obtain $\ \Gh^d:=\Gh\oplus\C\, d\ $
with Lie structure
$$[d\,,\,t]=0,\qquad [d\,,\,x(n)]=n\cdot x(n)\ .\tag\kn-3$$
The Virasoro algebra $\Lh$ is the Lie algebra with basis
$\ \{L_n\mid  n\in\Z\}\cup \{t_1\} $ and
Lie structure
$$[L_n\,,\,L_m]=(m-n)L_{m+n}+\delta^m_{-n}\frac {n^3-n}{12}\cdot
 \ce\cdot t_1,
\qquad [L_n\,,\,t_1]=0\ .\tag\kn-4$$
The number $\ce$ is sometimes called central charge.
We denote the Lie algebra without the element $t_1$ by $\L$.

Note that all these algebras are $\Z$-graded algebras if one
defines
$\deg(x(n))=n$, $\deg( L_n)=n$, and $\deg(t)=\deg(t_1)=\deg(d)=0$.
By the gradedness all above  algebras split
into subalgebras generated by the elements of negative,
 zero, resp\.  positive degree,
e.g\. $\ \Gh=\Gh_-\oplus \Gh_0\oplus\Gh_+$.
The finite dimensional algebra $\g$ can be embedded via
$\ x\mapsto x(0)$ into $\Gh$ and one obtains
$\ \Gh_0\cong \g\oplus\C\,t$.
This decomposition allows to define highest weight representations.
 These are
representations $V$ which are generated by one element
$\psi$ with $\G_+\psi=0$ and where the Cartan and the positive
nilpotent subalgebras of $\Gh_0$ operate in a certain manner.
In particular we get for every $v\in V$
and every $x\in\g$ that $x(n)v=0$ for $n$ big enough.
The similar construction works for the Virasoro algebra.

Indeed both algebras are related via the Sugawara construction. Let $V$
be a highest weight representation of $\Gh$, where
$\g$ is a simple Lie algebra,
$\{u_i,\;i=1,\ldots,\dim\g\}$ is a
basis and $\{u^i\}$ is the dual basis.
Let $\k$ be the dual Coxeter number (see Section \knsu\
for a description)
and suppose that $t$ operates as $\ce\cdot id$ on  $V$
and that  $\k+\ce\ne 0$ then
$$S_k:=-\frac 1{2(\k+c)}\sum_n\sum_i\nord{u_i(-n)u^i(n+k)}\tag\kn-5$$
is a well-defined operator. Here $u_i(n)$ is considered as operator
on $V$ and $\nord{....}$ denotes a normal ordering.
The normal ordering takes care that the elements of highest degree are
moved to the right where they eventually annihilate  every fixed vector.
What is quite astonishing is the fact that the
map $\ L_k\mapsto S_k\ $ and $\ t_1\mapsto id\ $ defines a
representation of the Virasoro algebra  with  central charge
$\ \dfrac {\ce\cdot\dim\g}{\ce+\k}\ $.

The highest weight of the Sugawara representation
we can be read off the relation
   $$-2(\ka+\ce)\cdot S_0\,\psi =(\l +2{\bar\rho}\,|\,\l )\,\psi\ ,
             \tag\kn-6
   $$
where $\l$ is the highest weight of  the $\Ghd$-module
from which we started
and
${\bar\rho}$ is the half-sum of the positive roots of the Lie algebra
$\fg .$

The above definitions look rather formal. But they have their geometric
interpretation.
The associative algebra $\ \C[z,z^{-1}]\ $ of Laurent polynomials can be
identified
with the algebra consisting of meromorphic functions on the Riemann sphere
$\widehat{\C}=\C\cup\{\infty\}$ (with quasi-global coordinate  $z\in\C$)
which are holomorphic outside $0$ and $\infty$.
In this picture $\G$ is identified with the space of
$\g$-valued meromorphic functions on $\widehat{\C}$
obeying the same regularity condition.
The algebra $\L$ can be interpreted as the  Lie algebra of
meromorphic vector fields on $\widehat{\C}$, again with this
regularity condition, if we identify $L_n$ with $\ z^{n+1}\fpz$.
In this description it is natural to ask for
generalization  to the  case of compact Riemann surfaces $\Gamma$ of
arbitrary  genus.
Essentially this generalization was done by Krichever and Novikov  in
\cite\KNFA. Their main reason was that in the usual quantization of
two-dimensional conformal field models by means
of the Virasoro algebra the role played by the
underlying Riemann surface is not clear.
They fixed two points $P_+$ and $P_-$ and considered
the algebra (resp\. Lie algebra) of meromorphic functions
(resp\. vector fields)
which are holomorphic outside the  two points.
(Indeed it is possible to consider the more  general algebra of
differential operators  \cite\SCHLCT.)
To obtain  central extensions they gave a geometric  definition
for the defining cocycles. Using this approach they were able to create
a far developed approach to the quantization on Riemann surfaces.

An essential step in the $g=0$ case was  to introduce  a
graded structure for these algebras.
This is not possible for higher genus. Fortunately in most cases
what is needed is a weaker concept, an almost-graded structure.
Krichever and Novikov did this by exhibiting  a certain basis, indexed
by the integers, and defining  the basis elements to be the homogeneous
elements. These basis elements will fulfil important duality relations.

We will describe this set-up in {\it Section 2}.
Generalized affine Kac-Moody algebras were
introduced by Krichever and Novikov \cite{\KNFA}, \cite {\KNJGP} and
extensively studied by  Sheinman \cite{\SHEIN}, \cite{\SHEINW},
\cite{\SHEINNS}
(see also Bremner \cite{\BREM},\cite\BREMa\ and Jaffe, Klimek and
 Lesniewski \cite{\JKL}).
For a multi-point generalization see Schlichenmaier \cite{\SCHLTH}, \cite
{\SCHLCT}.

As explained above one obtains in the classical case from
 such a representation by the
Sugawara construction a representation of  the Virasoro algebra
(with  central extension).
Our aim is to generalize this to higher genus with the goal to
obtain representations of a centrally extended
Krichever-Novikov vector field algebra.
This we will do in {\it Section 3},
where the main result is Theorem 3.1.
For $\g$  an abelian Lie algebra this has
been done by Krichever and Novikov \cite{\KNFA}.
In particular they made use of their beautiful technique of
delta-distributions on Riemann surfaces. In the
nonabelian case there is one
important point namely the appearance of the dual Coxeter number in the
final answer (see Theorem 3.1 below). It is necessary to  use earlier
ideas of \cite\KACRAI\  in order to explain this point in the framework of
the Krichever-Novikov approach.
For $\g$ a simple Lie algebra the  result has  been stated by
Bonora, Rinaldi, Russo and Wu in \cite{\BRRW}. There also a sketch of a
proof is given. The presentation there might not fulfil every
requirement of a scrupulous mathematician on a proof.  Especially if one
sees what delicate questions on normal ordering are involved.
  For the above
mentioned reasons and because in the appendix we want to generalize
the construction  to the multi-point situation  we present here a complete
proof (maybe we are following their line  of arguments).  By this  we
hope to convince also a  sceptical  mathematician that the result of
Bonora and collaborators is correct.
Again we obtain that the rescaled modes of the  Sugawara operator
of these representations will be a representation of a
certain (explicitly given) central extension of the vector field algebra
with the same central charge as in the genus zero case.
A key step in the proof, which is an important result by its own,
is
$$[S_k,x(n)]=x(\nabla_{e_k}A_n)\ ,\tag\kn-7$$
where $\nabla_{e_k}A_n$ is the function obtained by
applying the Lie derivative
with respect to the vector field $e_k$  to the function $A_n$,
where the index denotes the \KN\ degree
 (see Prop.~3.1 and Prop.~3.2).
This specializes for $g=0$ to
$\ [S_k,x(n)]=n\, x(n+k)$.

Some of the proofs will be postponed to {\it Section 4} which is in
 some sense
of more technical nature. An important tool in the proof will be the
duality property of  the Krichever-Novikov basis and the
''delta distribution'', see (2-19).
Nevertheless we prove there also very astonishing identities
between numbers which are obtained with the help of meromorphic
forms of different weights.
Up to now, we do not understand these identities completely.

In {\it Section 5} we recall  the notion of the weight of a
 representation of
the
Krichever-Novikov algebras of affine type
(generalized  affine Kac-Moody algebras) and of the Krichever
Novikov vector field algebra as introduced in  \cite\SHEINW .  We show
that with respect to the Sugawara  construction on  Riemann surfaces
the weights are related in a way  similar  to (1-6).
Additionally,
structure constants and the cocycle
of the algebra of meromorphic functions $\{A_n\}$
 are involved (see Theorem 5.1).

In {\it Section 6} we add a vector field $e$ to $\Gh $ that
 generalizes the
adjoining of the derivation $d$ in the classical case.  Let us denote
the obtained Lie algebra by $\Ghe$. We construct  higher genus Casimir
operators.  If the vector field $e$ is
given as a certain linear combination of the basis $e_k$  introduced by
Krichever and Novikov, then the Casimir operator $\Omega$ is given as $\
2L+2(\ce+\k)e$, where $L$ is a linear combination of the (not rescaled)
Sugawara operators $L_k$ with the same coefficients as in the
combination of $e$.
 Furthermore,
under the hypothesis that
the Casimir operator admits an eigenvector  we calculate the eigenvalue
in terms of the weight $\Lambda$
 calculated in  Section 5 and the highest weight of
the vector field $e$ in the corresponding $\Ghe$-module.

The Sugawara construction can be generalized to the
situation where one allows poles  at more than two points.
This will be done in {\it Appendix A}.
The crucial step for this more general set-up
is to introduce an almost-grading and
to find dual systems of basis elements.
This is done in \cite\SCHLTH, \cite{\SCHLL, 3.ref.}.
For a quick review see \cite\SCHLCT. (In this context see also
Sadov \cite{\SADOV}
and the appendix of \cite\KNPMS.)

\bigskip
The authors gratefully acknowledge the Volkswagen-Stiftung
for their support in the Research-in-Pairs program.
We would also like
to thank the Mathematisches Forschungsinstitut in Oberwolfach for
its hospitality while preparing this article.

%
\vskip 1cm
\head
2. The General Set-Up
\endhead
\def\kn{2}
\vskip 0.2cm
In this section we would like to recall the necessary facts on
the global Krichever-Novikov approach to conformal
field theory.
For the following let $\Gamma$ be a compact Riemann surface of
arbitrary genus $g$, $P_+$ and $P_-$ two fixed points which
for
genus $g\ge 1$ are in general position. Let $\Gamma^*=\Gamma\setminus
\{P_+,P_-\}$, and let $\rho$ be the unique meromorphic differential
with exact pole order 1 at the points $P_{\pm}$ and residues
$\res_{P_\pm}(\rho)=\pm 1$,
 holomorphic elsewhere and with purely imaginary periods.
We fix a point  $Q\in \Gamma^*$. The function
$\ u(P)=\re \int_Q^P\rho\ $
 is a well-defined harmonic function.
The level lines
$$C_\t=\{P\in\Gamma^*\mid u(P)=\t\}, \quad \tau\in\R\tag \kn-1$$
define a fibering of $\Gamma^*$.
For $\ \tau\ll 0\ (\tau\gg 0)\ $ the level line $C_\tau$ is a deformed
cycle around $P_+$, (resp\. $P_-$).

Let $K$ be the canonical bundle, i.e\. the bundle whose
local sections are the local holomorphic differentials. For every
$\l\in\Z$ we consider the bundle $K^\l:=K^{\otimes\l}$,
the
bundle with local sections the forms of
weight $\l$.
(After fixing a square root of the canonical bundle, a so-called
theta characteristic, it is possible to deal with $\l\in \frac 12\Z$.)
We denote
by $\Fl$ the vector space of global meromorphic sections of $K^\l$
 which are
holomorphic on $\Gamma^*$.
Special cases are  the differentials ($\l=1$),  the functions ($\l=0$),
and  the vector fields ($\l=-1$).
To denote  the space of functions we use also $\A$,
for the space of vector fields we  use also $\L$.

The (associative) algebra of functions $\A$
 operates by multiplication on $\Fl$.
The vector fields (i.e\. the elements in $\L$)
operate by taking the Lie derivative on $\Fl$.
In local coordinates
the Lie derivative can be described as
$$
\nabla_e(g)_|=(e(z)\fpz)\ldot (g(z)\dzl)=
\left( e(z)\pfz g(z)+\l\, g(z)\pfz e(z)\right)\dzl \ .\tag \kn-2
$$
Here and in the following we will use the same symbol for the
section of the bundle and its local representing function.
By (\kn-2) $\L$ becomes a Lie algebra and the vector spaces
$\Fl$ become Lie modules over $\L$ (i.e. we have
$ [\nabla_e,\nabla_f]=\nabla_{[e,f]}$).

Note that in the case $g=0$  with quasi-global coordinate $z$ and
$P_+=\{z=0\}$ and $P_-=\{z=\infty\}$ we obtain
$\A=\C[z,z^{-1}]$, the algebra of Laurent polynomials and for $\L$ the
Witt-Algebra,  the Lie algebra with basis
$\ \{\;l_n=z^{n+1}\fpz\mid n\in\Z\;\}$ and the commutator relation
$\ [l_n,l_m]=(m-n)\,l_{m+n}$.
We will call this case the ``the classical case''.

For the following let $\g$ be a finite dimensional reductive
Lie algebra, i.e\. a direct sum of an abelian and a semi-simple Lie
algebra, and let $\ (..|..)\ $ be a
 non-degenerate symmetric invariant bilinear form
on $\g$. In particular we have $([x,y]|z)=(x|[y,z])$.
We can take the Cartan-Killing form in the semi-simple
case and any non-degenerate symmetric form  in the abelian case.
The algebra $\G=\g\otimes \A$ is called the {\it current algebra}.
It can be considered as the algebra of $\g$-valued meromorphic
functions on
$\Gamma$ which are holomorphic on $\Gamma^*$.

We take an element $\a\in \HG$ and represent it by a cycle
$Z_\a$ which is a sum of differentiable curves on $\Gamma$.
The map
$$\gamma_\a:\A\times\A\to\C,\qquad
\gamma_\a(f,g):=\cintb fdg\tag\kn-3$$
defines a 2-cocycle for the abelian Lie algebra $\A$.
Due to the fact that all poles are located in $P_\pm$ homologous cycles
$Z_\a$ define the same cocycle $\gamma_\a$.
Using the cocycle (\kn-3) a central extension $\Ah_\a$ is obtained:
$\
0\ \mapright {}\ \C \ \mapright {}\
\Ah_\a\ \mapright {}\ \A\ \mapright{} \ 0$.
If $\hat f$ and $\hat g$ are lifts of elements $f,g\in\A$ and $t$ is
a generator of the center then
$$[\hat f,\hat g]=-\gamma_\a(f,g)\cdot t\ ,\qquad
[\,t,\Ah_a]=0\ .\tag \kn-4$$
The algebra $\Ah_\a$ is called a {\it Heisenberg algebra (of higher genus)}.
This can be generalized to the current algebra. We take
$\Gh_\a=\G\oplus \C\cdot t$ as vector space and define the Lie algebra
structure as follows. For $a\in\G$ we use $\hat a=(a,0) $ for the
corresponding lift in $\Gh_a$ and define for $x,y\in\g$
$$[\widehat{x\otimes f},\widehat{y\otimes g}]=
\widehat{[x,y]\otimes (f g)}-(x|y)\cdot\gamma_\a (f,g)\cdot t,\qquad
[\,t,\Gh_\a]=0\ .\tag\kn-5$$
{}From the invariance of the bilinear form $(..|..)$ it follows that
this is indeed a central extension of $\G$.
The algebras obtained in this way
are called {\it \KN\ algebras of affine type}. They
are higher genus generalizations of
affine Kac-Moody algebra.
In Section 6 we will adjoin an additional vector field
to the algebra $\Gh_\a$.

In the following an  important role will be played by the central
extensions  corresponding to the cycle which is represented by
a level line  $C_\tau$. If we use $\gamma_0$ or just $\gamma$ we will
always mean this 2-cocycle.
And we will denote by $\Ah$ and $\Gh$ the corresponding central
extensions.

In this article we also have to deal  with central extensions of the
vector field  algebra. Krichever and Novikov  gave a generalization
of the standard cocycle of the Virasoro algebra to higher genus as follows.
Let $R$ be a  projective connection (holomorphic or meromorphic
with poles only at $P_\pm$)
on $\Gamma$.
For the definition of a projective connection see \cite\KNFA.
Note that the difference of two projective connections is
a form of weight 2 (i.e\. a quadratic differential).

 For vector
fields $e$ and  $f$ represented locally as
$e(z)\fpz$ and $f(z)\fpz$ the 2-cocycle is defined as
$$\chi_{\a,R}(e,f)=\cinttb \left(
\frac 12(e'''f-ef''')-R\cdot(e'f-ef')\right) dz\ .
\tag \kn-6$$
Without the connection $R$ the expression under the integral
would not be a well-defined differential.
The choice of a different projective
connection will
yield a cohomologous cocycle.
The central extension is  given  by  $\Lh_\a=\L\oplus \C\cdot t$
as vector space  with Lie structure (using $\hat e=(e,0)$)
$$[\,\widehat{e},\widehat{f}\,]=\widehat{[e,f]}+\chi_\a(e,f)\cdot t,\qquad
[\,t,\Lh]=0\ .\tag\kn-7$$
Again we use $\Lh=\Lh_0=\Lh_{[C_\tau]}$ if we integrate over a
level line.
There is also a suitable extension to the multi-point situation
(see Appendix A, \cite{\SCHLTH} and the third article in Ref.
 \cite{\SCHLL}).

For the construction of  highest weight representations in the
classical  case it is important that the above  algebras are graded
 algebras.
Note that in this case there is only one cycle class, hence only
one  nontrivial central extension (up to equivalence
and isomorphy) defined as above. In the higher genus  it is the concept
of almost-grading which will do the job.
Such a grading has been introduced by Krichever  and Novikov in the
following way.
Let $g=0$  and $\l\in\Z$ be arbitrary or $g\ge 2$  and $\l\ne 0,1$ then
there is for every $n\in\Z$ a unique (up to multiplication with a
scalar)  $f_n^\l\in \Fl$ such that
$$\ord_{P_+}(f_n^\l)= n-\l,\qquad
\ord_{P_-}(f_n^\l)= -n+(\l-1)+(2\l-1)(g-1)\ .\tag\kn-8$$
If we fix a local coordinate $z_+$ at $P_+$  we adjust the scalar
by requiring  locally
$$f_n^\l(z_+)_|=z_+^{n-\l}(1+O(z_+))\dzl\ .\tag\kn-9$$
Note that our indexing differs from the one used by
Krichever and Novikov by a shift.
The set $\{f_n^\l\mid n\in\Z\}$ is a basis of $\Fl$.
By calculation of residues at $P_\pm$ we obtain
$$\cintt f_n^\l\cdot f_{m}^{1-\l}=\delta_{n,-m}\ .\tag\kn-10$$
It is quite convenient to introduce the
notation
$\ f_\l^{*,n}=f_{-n}^\l,\
A_n=f_n^0,\  e_n=f_n^{-1},\ \omega^n=f_1^{*,n},
\ \Omega^n=f_2^{*,n}\ $.
Now the duality reads as
$$\cintt A_n\omega^m=\delta_n^m,\qquad
\cintt e_n\Omega^m=\delta_n^m\ .\tag\kn-11$$
For the remaining cases of $\l$ and $n$ we have to modify the
above  prescription
for $-g \le n\le 0$.
We set $A_0=1$ and $\w^0=\rho$ and fix the elements
$A_n$ and $\w^n$ for $-g\le n\le -1$  by
$$\alignedat 2
\ord_{P_+}(A_n)&=n,&\qquad
\ord_{P_-}(A_n)&=-n-g-1\ ,\\
\ord_{P_+}(w^n)&=-n-1,&\qquad
\ord_{P_-}(w^n)&=n+g\ ,
\endalignedat\tag \kn-12$$
 and the duality (\kn-11).
In all cases we obtain the same order at the point $P_+$ and the
same duality relation as in the generic case.
We will use the term {\it critical strip} to denote the
index values $\ -g,\ldots,-1,0$.

We define $\deg f_n^\l:=n$  and call this elements homogeneous
elements of degree  $n$. With respect to this degree the
algebras $\ \A,\ \L,\ \G\ $ are almost-graded (where in the latter case
we define  for $x\in\g,\  \deg (x\otimes A_n):=n$)
and the vector spaces $\Fl$ are almost-graded modules over $\A$ and $\L$.
More precisely, we have
$$A_n\cdot A_m=\sum_{k=n+m}^{n+m+L}\a_{nm}^k\,A_k,\qquad
[e_n,e_m]=\sum_{k=n+m}^{n+m+M} C_{nm}^k\,e_k\ ,\tag\kn-13$$
where the constants are given by the duality relations as
$$
\a_{nm}^k=\cintt A_nA_m\w^k,\qquad
C_{nm}^k=\cintt ([e_n,e_m])\cdot\Omega^k\ ,\tag\kn-14$$
and the constants $L$ and $M$ do not depend on $n$ and $m$.
For $g\ne 1$ we calculate  $M=3g$.
For $g=0$ we obtain $L=0$. In general
explicit formulas can be given. For us of importance is only that
 if $n$ and $m$ are both
on  the same  side of the critical strip
(e.g\. $n,m<-g$) then the upper bound will be $n+m+g$.

We want to extend our almost-grading to the central
extensions $\ \Ah_\a,\,\Lh_\a,\,\Gh_\a\ $ by defining
$\deg\widehat{x}:=\deg x$ and $\deg t:=0$.
For this to work  our 2-cocycles which define the central
extensions should be  local cocycles \cite{\KNFA}, i\. e\. there
should  be constants $K$ and $N$, such that
$$\gamma_\a(A_n,A_m)=0,\quad\text{for}\ |n+m|>K,\qquad
\chi_\a(e_n,e_m)=0,\quad\text{for}\ |n+m|>N\ .\tag\kn-15$$
For arbitrary $\a\in\HG$ this will not be the case. But if we integrate
over a level line $C_\tau$
we get
that $\gamma(A_n,A_m)\ne 0$ implies $\ -2g-2\le n+m\le 0$
(for generic $n$ and $m$ we get as lower bound even $-2g$) and that
 $\chi(e_n,e_m)\ne 0$ implies $\ -6g\le n+m\le 0$
(for $g\ne 1$).
For later reference we note that at the upper boundary we
obtain
$$\gamma(A_n,A_{-n})=(-n),\qquad
\chi(e_n,e_{-n})=\frac 1{12}(n^3-n)\ .\tag\kn-16$$

We consider the vector space decomposition
$$
\gathered
\A=\A_-\oplus\A_0\oplus\A_+,\qquad\text{with}\\
\A_-:=\langle A_n\mid n\le -g-1\rangle\ ,\quad
\A_0:=\langle A_n\mid -g\le  n\le 0\rangle,\quad
\A_+:=\langle A_n\mid n\ge 1\rangle,
\endgathered\tag\kn-17$$
and the corresponding decomposition
$$\G=\G_-\oplus\G_0\oplus\G_+,\quad\text{with}\quad
\G_\beta=\g\otimes \A_\beta,\quad \beta\in\{-,0,+\}\ .\tag\kn-18$$
Due to the almost-gradedness of $\A$ we see that $\A_+$ and $\A_-$, resp\.
$\G_+$ and $\G_-$ are subalgebras. Contrary to the classical case
for higher  genus $\A_0$ and $\G_0$ will only be subspaces.
We will call the elements of $\A_0$ also the elements from the
critical strip.
By the locality of the cocycle this decomposition extends to
$\ \Gh=\Gh_-\oplus\Gh_0\oplus\Gh_+\ $, where
$\Gh_\pm$ can be identified with $\G_\pm$ and $\Gh_0=\G_0\oplus\C\cdot t$.
Clearly $\Gh$ is generated by $\widehat{x\otimes A_n}$ for $x\in\g$ and
$n\in\Z$ and by the central element $t$.
We will denote
the element $\widehat{x\otimes A_n}\in\Gh$ and  ${x\otimes A_n}\in \G$
also by $x(n)$ if convenient. Note that for $g=0$ we have
$A_n=z^n$ and
our notation specializes completely to the classical notation.

We close this section by introducing the very useful object
$$\Delta(Q',Q)=\sum_{n\in\Z}A_n(Q')\w^n(Q)\ .\tag\kn-19$$
It can be considered  as the {\it delta distribution} in the sense
that we have for $f\in\A$ and $\w\in\Fn 1$
$$\cintt \Delta(Q',Q)f(Q)=f(Q'),\quad
\text{resp\.}\quad
\cinttp \Delta(Q',Q)\w(Q')=\w(Q)\ .\tag\kn-20$$
Of course (\kn-19) can be extended to arbitrary pairs of weights
$\ (\l,1-\l)\ $.
%
%
\vskip 1cm
\head
3. The Sugawara Construction
\endhead
\def\kn{3}
\vskip 0.2cm
\definition{Definition}
A module $V$ over the Lie algebra $\Gh$
(resp\. a representation) is called an {\it admissible
module} (resp\. {\it representation}) if for every $v\in V$
and for all $x\in\g$ we have $x(n)v=0$ for $n\gg 0$.
\enddefinition
Let $V$ be a fixed admissible module.
Let the central element $t$ operate by multiplication with a
scalar $\ce\in\C$.
In the classical case these are the usual highest weight modules
of affine Kac-Moody algebras. For higher genus such modules have been
studied by Sheinman \cite{\SHEIN}, \cite{\SHEINW}, \cite{\SHEINNS}.
If $x\otimes A$ (or more precisely $\widehat{x\otimes A}$)
is an element of $\Gh$ then we will use the notation $\ x(A)\ $
for the corresponding operator on $V$. For
$x\otimes A_n$ with $A_n$ the special basis elements we will
also use for short $\ x(n)\ $ to denote $x(A_n)$.

Recall that we assume $\g$ to be a finite dimensional
reductive Lie algebra. We choose a
basis $\ u_i,\ i=1,\ldots,\dim\g\ $ of $\g$ and the corresponding
dual basis $\ u^i,\  i=1,\ldots,\dim\g$ with respect to the invariant
non-degenerate  symmetric bilinear form $(..|..)$.
The Casimir  element $\ \Omega^0=\sum_{i=1}^{\dim \g} u_iu^i\ $
of the universal enveloping algebra  $U(\g)$ is independent
of the choice of the basis.
In the following a summation over $i$ is always assumed to be
over the above summation range.
\proclaim{Lemma 3.1}

\noindent
(1)\qquad $[\,\Omega^0,\g]=0$ .

\noindent
(2)\qquad $\sum_{i}[u_i,u^i]=0$ .

\noindent
(3)\qquad $\sum_i[u_i\otimes A_n,u^i\otimes A_m]=
-\dim\g\cdot\gamma(A_n,A_m)\cdot  t$ .

\noindent
(4)\quad For  $\g$ an abelian or a simple Lie algebra  there is a
constant $\k$,
such that \nl
$\ \sum_i ad_{u_i}\circ ad_{u^i}=2\k$.
\endproclaim
\demo{Proof}
(1) and (4) are standard knowledge. (2) follows from the invariance
of bases.
To show (3) we take the structure equation of the centrally extended
algebra
$$[u_i\otimes A_n,u^i\otimes A_m]=
[u_i,u^i]\otimes (A_nA_m)-\gamma(A_n,A_m)\cdot
(u_i|u^i)\cdot t\ .$$
After summation over $i$ the first summand will vanish by (2) and we get
the result.
\qed
\enddemo
Note that $2\k$ is the eigenvalue of the Casimir operator in the
adjoint representation. In the
case where $\g$ is simple $\k$ is the {\it dual Coxeter number}.
In the abelian case $\ \k=0$.
Let us define for $Q\in\Gamma$ the formal sum (the ``generating
function'')
$$\widehat{x}(Q)=\sum_{n}x(n)\cdot \w^n(Q)\ .\tag\kn-1$$
Here and in the following a summation over the indices of our
forms will always mean a summation over $\Z$ if not stated
otherwise.
We define the higher genus {\it Sugawara} (or {\it Segal}) {\it operator}
$$T(Q):=\frac 12\sum_i\nord{\widehat{u_i}(Q)\widehat{u^i}(Q)}\quad =\
\frac 12\sum_{n,m}\sum_i\nord{u_i(n)u^i(m)}\w^n(Q)\w^m(Q)\ .\tag\kn-2$$
Here $\ \nord{....}\ $ denotes some normal ordering.
If we consider again $T(Q)$ as ``generating function''
we can write
$$T(Q)=\sum_kL_k\cdot\Omega^k(Q)\tag\kn-3$$
with certain operators $L_k$. Using duality we obtain
$$
\gathered
L_k=\cintt T(Q)e_k(Q)=\frac 12\sum_{n,m}\sum_i
\nord{u_i(n)u^i(m)}l_k^{nm},\\
\text{with}\qquad
\ l_k^{nm}=\cintt w^n(Q)w^m(Q)e_k(Q)\ .
\endgathered\tag\kn-4$$
Note that for a fixed value of $k$ for every value  of $n$ there is
only a finite set of values for $m$ such that $l_k^{nm}\ne 0$.
More precisely,  $l_k^{nm}\ne 0 $ implies that for the
indices  $n$ and $m$ outside the
exceptional strip  we have  $k-n\le m\le k-n+g$. This specializes
 for $g=0$ to
$l_k^{nm}=\delta_k^{m+n}$  which  gives the usual definition of the
$g=0$ Sugawara  operators (see \cite\KACRAI\  and references therein).
By this finiteness and the normal ordering the operators $L_k\in gl(V)$
 are well-defined.

 We will choose in Section 3 and Section 4
the prescription
$$
\nord{x(n)y(m)}\quad :=\ \cases x(n)y(m)&,n\le m\\
                           y(m)x(n)&,n>m\
                     \endcases\tag\kn-5
$$
\noindent as normal ordering.
We will see in the proofs that
Proposition~\kn.1 and hence Proposition~\kn.2 will not depend
on the normal ordering. In Proposition~\kn.3 where we show that
the Sugawara operators define a centrally extended
Krichever-Novikov algebra we will see that only
the cohomology class of the cocycle defining the
central extension will depend on the
normal ordering  chosen.

The following proposition is the key step in the construction.
\proclaim{Proposition 3.1}
Let $\g$ be either an abelian or a simple Lie algebra, then
$$\gather
[L_k,x(r)]=-(\ce+\k)\sum_vK_{r,k}^vx(v)\ ,\tag\kn-6\\
\text{with}\quad
K_{r,k}^v:=\cintt w^v(Q)e_k(Q)dA_r(Q)=
\sum_m l_k^{vm}\gamma_{mr},\tag\kn-7\\\
\gamma_{mr}:=\gamma(A_m,A_r):=\cintt  A_m(Q)dA_r(Q)\ .
\tag\kn-8
\endgather
$$
The result does not depend on the normal ordering.
\endproclaim
\noindent
Note that all the infinite sums above are indeed well-defined finite
sums. The above proposition will specialise
in the classical case to
$\ [L_k,x(r)]=-(\ce+\k)r\,x(r+k)$ (see \cite{\KACRAI, Prop.10.1}).
We will postpone the proof to the next section. Here we want to
show the relation (\kn-7). Using the ``delta distribution'' (2-19)
 we obtain
$$\gather
\sum_m l_k^{vm}\gamma_{mr}=\sum_m\cintt \w^v(Q)\w^m(Q)e_k(Q)
\cint {\tau'}A_m(Q')dA_r(Q')=
\\
\cintd \w^v(Q)e_k(Q)dA_r(Q')\Delta(Q',Q)=
\cint {\tau'}\w^v(Q')e_k(Q')dA_r(Q')=K_{r,k}^v\ .\qed
\endgather$$
We define the operation of $e\in\L$ on $\widehat{x}(Q)$ as
$$e\ldot \widehat{x}(Q):=\sum_n x(n) (\nabla_{e}\w^n)(Q),\qquad
\text{for}\quad
\widehat{x}(Q)=\sum_nx(n)\w^n(Q)\ .\tag\kn-9$$
\proclaim{Proposition  3.2}

\noindent
(1)\qquad $[L_k,x(r)]=-(\ce+\k)\,x(\nabla_{e_k}A_r)$ .

\noindent
(2)\qquad  $[L_k,\widehat{x}(Q)]=(\ce+\k)e_k\ldot \widehat{x}(Q)$ .
\endproclaim
\demo{Proof}
We can write $ \nabla_{e_k}A_r$ in local coordinates as
$e_k(z)\pfz {A_r(z)}$, hence
$$ (\nabla_{e_k}A_r)(Q)=e_k(Q)dA_r(Q)=\sum_v\beta_{r,k}^vA_v(Q)\ .$$
By duality the coefficients $\beta_{r,k}^v$ calculate as
$\ \beta_{r,k}^v= \cintt e_k(Q)dA_r(Q)\w^v(Q)=K_{r,k}^v$.
Hence Prop. 3.1 implies (1).
\nl
To prove (2) we write
$\ \nabla_{e_k}\w^v=\sum_r\zeta_{r,k}^v\w^r\ $ where
$$\zeta_{r,k}^v=\cintt(\nabla_{e_k}\w^v)A_r=
\cintt\nabla_{e_k}(\w^vA_r)-\cintt\w^v\nabla_{e_k}A_r=-K_{r,k}^v\ ,$$
because the residue of a Lie derivative of a  meromorphic
differential will vanish (see for example \cite{\SCHLTH, p.102}).
This proves (2). \qed
\enddemo
\proclaim{Proposition 3.3}
The operators $L_k\in gl(V)$ and $id=1\in gl(V)$ close up to
a Lie subalgebra of $gl(V)$ with commutator relation
$$[L_k,L_{l}]=-(\ce+\k)\sum_nC_{kl}^nL_n-\frac 12
\ce(\ce+\k)\dim \g\cdot \chi_{kl}\cdot id
\tag\kn-10$$
where  $C_{kl}^n$ are the structure constants of the vector field
algebra  $\L$ and
$$\gathered\chi_{kl}= \psi_{kl}+\widehat{\chi}_{kl},\quad
\psi_{kl}=\sum_{s,v}\sum_{n=0}^{v+1} C_{kl}^s l_s^{nv}\gamma_{nv}
=\sum_{v}\sum_{n=0}^{v+1} E_{kl}^{nv}\gamma_{nv},\\
E_{kl}^{nv}=\cintt [e_k,e_l]\cdot\w^n\w^v,\qquad
\widehat{\chi}_{kl}=\bigg(\sum\Sb n>0 \\ v\le 0 \endSb
-\sum\Sb n\le 0 \\ v> 0 \endSb
\bigg) K_{v,k}^n K_{n,l}^v
\ .
\endgathered\tag\kn-11$$
The $\chi_{kl}$ could be non-zero only if
$\ -6g\le k+l\le 0$.
A different normal ordering will change the range of the
$n$ and $v$ summation in the definition of  $\psi_{kl}$
and hence change $\chi_{kl}$.
At the upper bound we obtain
$$\chi_{k,-k}=-\frac 16(k^3-k)\ .\tag\kn-12$$
\endproclaim
\noindent If $\ \ce+\k\ne 0\ $ we can choose the rescaled  elements
$\ L_k^*=\frac {-1}{\ce+\k}L_k\ $ and obtain
$$[L_k^*,L_l^*]
=\sum_n C_{kl}^n L_n^*- \frac {\ce}{2(\ce+\k)}\dim\g
 \cdot \chi_{kl}\cdot id\
.\tag\kn-13$$
Hence we obtain a representation of a  centrally extended algebra
of $\L$. Note that from the Jacobi identity inside $gl(V)$ and inside
$\L$ it follows
that $\chi_{kl}$ indeed defines a 2-cocycle $\chi$ for the
 Lie algebra $\L$.
By Proposition~3.3 this cocycle is local. Krichever and Novikov showed
(\cite{\KNFA}, \cite{\KNJGP}) that all local cocycles are cohomologous
to a multiple of the geometric cocycle $\chi_R$ (2-6) with
a suitable projective connection
$R$ (and integration over a level line $C_\tau$).

Hence the centrally extended algebra is a representation  of the
algebra $\Lh$.
To study it in more detail the following facts are quite useful.
\proclaim{Lemma~\kn.2}
Let $\chi_R$ be the cocycle (\knset-6) (where the integration curve
equals $C_\tau$).
\vskip 0.1cm
\noindent (a) $R$ is a meromorphic projective connection with only poles
up to order two at the points $P_\pm$ if and only if
$\chi_R(e_k,e_{l})=0$ for $k+l>0$ or  $k+l <-6g$.
\vskip 0.1cm
\noindent (b)
Let $R_|(z_+)=\alpha_+z_+^{-2}(1+O(z_+))$ be the local form of the
projective connection at the point $P_+$, then
$$\chi_R(e_k,e_{-k})=\frac 1{12}(k^3-k-2\a_+ k)$$
\endproclaim
\noindent
To show this one calculates the involved residues (see also \cite\KNFA).

\noindent
By Proposition~\kn.3 we see that $\chi$ is a cocycle which fulfils the
conditions of Lemma~3.2(a)
Hence $\ \chi=d\cdot\chi_R$ with $d\in\C$ and $R$ a suitable
projective connection. If we compare
$\ \chi_R(e_k,e_{-k})\ $ and $\ \chi(e_k,e_{-k})\ $ we see that
$\a_+=0$ (hence $\ord_{P_+}(R)\ge -1$) and $d=-2$.
Altogether we obtain the following theorem:
\proclaim{Theorem 3.1}
Let $\g$ be either a finite dimensional abelian or simple Lie algebra and
 $2\k$ be the eigenvalue of the Casimir operator in the adjoint
representation and  $\Gh$ be the higher genus affine Kac-Moody algebra.
Let $V$ be  an admissible representation  where the central element
operates as  $\ce\cdot identity$. If $\ce+\k\ne 0$ then the rescaled
modes
$$L_k^*=\frac {-1}{2(\ce+\k)}\sum_{n,m}\sum_i \nord{u_i(n)u^i(m)}l_k^{nm}
=\frac {-1}{2\pi\i(\ce+\k)}\int_{C_\tau} T(Q) e_k(Q)\ ,$$
of the Sugawara operator define a representation of
a central extension of the Krichever-Novikov vector field
algebra given by
the geometric cocycle
$$\chi(e,f)=\frac {\ce\cdot\dim\g}{(\ce+\k)}\cdot
\cinttt \left(
\frac 12(e'''f-ef''')-R\cdot(e'f-ef')\right) dz\ ,
\tag\kn-15
$$
with a suitable meromorphic projective connection $R$ with poles only at
$P_\pm$ and
\nl
$\ord_{P_+}(R)\ge -1$ and $\ord_{P_-}(R)\ge -2$.
\endproclaim
\remark{Remark}
The expressions $\psi_{kl}$ in (\kn-11) define  just a coboundary
(in the sense of Lie algebra cohomology). Recall a 2-cocycle $\psi$
is a coboundary of the Lie algebra  $\L$ if there is a linear form
$\Phi:\L\to\C$ such that   $\ \psi(e,f)=\Phi([e,f])$.
We define
$\Phi$ by $\ \Phi(e_s):=\sum_v\sum^{v+1}_{n=0} l_s^{nv} \gamma_{nv}\ $
(this is a finite sum)
and calculate
$\ \Phi([e_k,e_l])=\Phi(\sum_s C_{k,l}^s e_s)=\psi_{kl}$.
Again a different normal ordering would result in a different
range of summation in the definition of $\Phi$ above.
By this we proved again, without using the result of Krichever and
Novikov on the local cocycles, that different normal ordering
would not change the cohomology class of the central extension.
\endremark
%
%
\vskip 1cm
\head
4. The Proof of Proposition 3.1 and  Proposition 3.3
\endhead
\def\kn{4}
\vskip 0.2cm
In the definition of the $L_k$ formal infinite sums
of operators are involved. To take care about the well-definedness  we use
cut-off functions as has been done by Kac and Raina in \cite{\KACRAI}.
Let $\psi $ be the function on $\R$ given as
$$\psi(x)=1\quad \text{if}\ |x|\le 1\qquad\text{and}\qquad
\psi(x)=0\quad \text{if}\ |x| >1\ .\tag\kn-1$$
For $\ep\in \R$ we define
$$L_k(\ep)=\frac 12\sum_{n,m}\sum_i\nord{u_i(n)u^i(m)}
l_k^{nm}\psi(\ep n)\ .
\tag\kn-2$$
We fix $k$.
For every $n$ there are only finitely many $m$ such that $l_k^{nm}\ne 0$.
Hence for $\ep>0$ the sum consists only of finitely many summands.
If $v\in V$ then by the normal ordering prescription only finitely many
operators
$l_k^{nm}\nord{u_i(n)u^i(m)}$ will operate non-trivially  on $v$.
Hence if we choose $\ep>0$ small enough we get
$L_k(\ep)v=L_kv$.
This we will mean if we write  $\lim_{\ep\to 0}L_k(\ep)=L_k$.

If we drop the normal ordering symbols in (\kn-2)  we obtain
an expression $\widetilde{L}_k(\ep)$  which is well-defined as long
as $\ep\ne 0$.
For every pair $(n,m)$ which is not in normal order we take up
the commutator
$\sum_i[u_i(n),u_i(m)]$ which is a scalar by Lemma~3.1(3), hence
$\ L_k(\ep)=\widetilde{L}_k(\ep)+\alpha\cdot t$,
where $\alpha$ is a scalar as  long as $\ep\ne 0$.
In particular,
if we calculate  commutators we can forget about the normal ordering
as long as we stay with $\ep\ne 0$.

\subheading{Proof of Proposition 3.1}
Per definition we have
$$\multline R_\ep:=2\,[\widetilde{L}{}_k(\ep),x(r)\,]=
\sum_{n,m}\sum_i [u_i(n)u^i(m),x(r)]\, l_k^{nm}\psi(\ep n)
\\
=\sum_{n,m}\sum_i \left(u_i(n)[u^i(m),x(r)]+[u_i(n),x(r)]u^i(m)\right)
 l_k^{nm}\psi(\ep n),
\endmultline$$
after expanding the commutator and reordering the elements
again. Each commutator can now be written like
$$[u^i(m),x(r)]=[u^i,x](A_mA_r)-(u^i|x)\gamma_{mr}\cdot \ce$$
(note that $t.v=\ce\cdot v$).
Hence we obtain
$\ R_\ep=A_\ep+B_\ep-(C_\ep+D_\ep)\ $ where
$$\alignat 2
A_\ep&=\sum_{n,m}\sum_i u_i(n)[u^i,x](A_mA_r) l_k^{nm}\psi(\ep n),
\
&B_\ep&=\sum_{n,m}\sum_i [u_i,x](A_nA_r)u^i(m)l_k^{nm}\psi(\ep n),
\\
C_\ep&=\sum_{n,m}\sum_i u_i(n)(u^i|x)\gamma_{mr}l_k^{nm}\ce\psi(\ep n),\
&D_\ep&=\sum_{n,m}\sum_i (u_i|x)u^i(m)\gamma_{nr}l_k^{nm}\ce\psi(\ep n)\ .
\endalignat
$$
Now using $\ \sum_i u_i\otimes A_n(u^i|x)= \left(\sum_i (u^i|x)u_i\right)
\otimes A_n=
x\otimes A_n=x(n)\ $,
we obtain
$\
C_\ep=\sum_{n,m}x(n)\gamma_{mr}l_k^{nm}\ce\psi(\ep n)\ $, and
$\
D_\ep=\sum_{n,m}x(m)\gamma_{nr}l_k^{nm}\ce\psi(\ep n)
\ $.
For fixed $r$ and $k$ only finitely many terms  occur. Hence
for
$\ep=0$
we obtain
$$\lim_{\ep\to 0} (C_\ep+D_\ep)=
2\,\ce\cdot\sum_n\big(\sum_m l_k^{nm}\gamma_{mr}\big)x(n)=
2\,\ce\sum_n K_{r,k}^n x(n)\ .\tag\kn-3$$
Here we used (\knsu-7).
The sums $A_\ep$ and $B_\ep$ for $\ep\to 0$ do not make sense separately.
We can make sense out of them if we change to normal ordering.
For this we have to resolve
$A_mA_r=\sum_s \a_{mr}^sA_s$ with $\ \a_{mr}^s=\cintt A_mA_r\w^s$.
We obtain
$$A_\ep=\sum_{n,m,s}\sum_i u_i(n)[u^i,x](s) \a_{mr}^s l_k^{nm}\psi(\ep n),
\quad
B_\ep=\sum_{n,m,s}\sum_i [u_i,x](s)u^i(m)\a_{nr}^s l_k^{nm}\psi(\ep n)
\ . $$
For the elements which are not in normal order we have to
pick up a commutator. We write $A_\ep=A_\ep^{(1)}+A_\ep^{(2)}$ and
$B_\ep=B_\ep^{(1)}+B_\ep^{(2)}$
where  $A_\ep^{(1)}$ resp\. $B_\ep^{(1)}$ are the expressions above
just with the normal ordering columns.
We can write the commutator
$$[u_i(n),[u^i,x](s)]=[u_i,[u^i,x]](A_nA_s)-\gamma_{ns}(u_i|[u^i,x])\ce\ .$$
If we sum over $i$ the second term will vanish because
$\ (u_i|[u^i,x])=([u_i,u^i]|x)\ $ and Lemma 3.1(2).
Lemma 3.1(4) gives for the first term
$\ 2\k\cdot x(A_nA_s)=2\k \sum_v \a_{ns}^v x(v)$.
Applying the same to $B_\ep^{(2)}$ we get
$$
A_\ep^{(2)}+B_\ep^{(2)}= 2\k\sum_v \bigg(\sum_{s,m}\sum_{n>s}
\a_{ns}^v \a_{mr}^s l_k^{nm}\pem-
\sum_{n,m}\sum_{s>m}
 \a_{sm}^v \a_{nr}^s l_k^{nm}\pem
\bigg) x(v)\ .\tag \kn-4$$
Note that neither sum alone will make sense if we put $\ep=0$.
Before we continue to deal with (\kn-4) we first show that
$A_0^{(1)}+B_0^{(1)}$ will vanish.
First we change  the variables in  the summation
for $B_\ep^{(1)}$ in the way $s\to n\to m\to s$.
By the normal ordering $A_0^{(1)}$ and $B_0^{(1)}$
are well-defined operators
in the sense that applied to a fixed $v\in V$ only finitely many
summands will act nontrivially. Hence we can forget about
the $\psi-$factor.
 \proclaim{Lemma \kn.1}
For $\ F_{r,k}^{sn}:=\sum_m\a_{mr}^sl_k^{nm}=\cintt
 A_r(Q)\w^s(Q)\w^n(Q) e_k(Q)\ $
we have
$\ F_{r,k}^{sn}=F_{r,k}^{ns}\ $.
\endproclaim
\demo{Proof}
$$\multline
F_{r,k}^{sn}=\sum_m\cintd A_m(Q)A_r(Q)\w^s(Q)\w^n(Q')\w^m(Q')e_k(Q')
\\=
\cintd  A_r(Q)\w^s(Q)\w^n(Q')e_k(Q')\Delta(Q,Q')=
\cintt A_r(Q)\w^s(Q)\w^n(Q) e_k(Q)\ .
\endmultline$$
This is obviously symmetric in $n$ and $s$.\qed
\enddemo
Now
$\ A_0^{(1)}+B_0^{(1)}=
\sum_{n,s}\sum_i\left(\nord{u_i(n)[u^i,x](s)+[u_i,x](n)u^i(s)}
\right)F_{r,k}^{sn}\ .$
\proclaim{Lemma \kn.2}
$\qquad
\sum_i\left(\nord{u_i(n)[u^i,x](s)+[u_i,x](n)u^i(s)}\right)=0\ $.
\endproclaim
\demo{Proof}
We calculate
$\
\sum_i u_i(n)[u^i,x](s)=\sum_i u_i(n)\sum_j([u^i,x]|u_j)u^j(s)=$\nl$
-\sum_{i,j} u_i(n)(u^i|[u_j,x])u^j(s)=
-\sum_{j} [u_j,x](n)u^j(s)\ $.\qed
\enddemo
\noindent Hence the $\ A_0^{(1)}+B_0^{(1)}=0$.

We now take up  (\kn-4) again.
\proclaim{Claim}
The expression inside the $\ v $-summation is for $\lim_{\ep\to 0}$ equal
to
the $\lim_{\ep\to 0}$ of
$$
E_{\ep}^{(N)}:=\sum_{m,s}\sum_{n>N}\a_{ns}^v \a_{mr}^s l_k^{nm}\pem-
\sum_{n,m}\sum_{s>N}\a_{sm}^v \a_{nr}^s l_k^{nm}\pem\ ,\tag\kn-5$$
where $N$ is an arbitrary  integer.
\endproclaim
\demo{Proof}
If we calculate the difference we obtain
\footnote{
If the upper bound is smaller than the lower bound, we mean that
one has to
switch the summation range and the sign of the expression.
This will be understood in all the summations which will follow.
}
$$
\sum_{m,s}\sum_{n=s+1}^N\a_{ns}^v \a_{mr}^s l_k^{nm}\pem-
\sum_{n,m}\sum_{s=m+1}^N\a_{sm}^v \a_{nr}^s l_k^{nm}\pem\ .\tag\kn-6$$
Note that due to the almost-grading
in each sum for fixed $v,k,r$ only finitely many terms are involved.
Hence we can forget about $\pem$ and change variables ($s\to n\to m\to s$)
in the second sum. Applying Lemma \kn.1
we see that the difference will vanish.
\qed\enddemo
\noindent This proof shows also  that our  result will not depend on
the normal ordering chosen. Again the difference will
consist of finitely
many terms which will cancel.

We will study
$$E_\ep^{(0)}=\sum_s\sum_{n>0}\a_{ns}^v F_{r,k}^{sn}\pem-
\sum_n\sum_{s>0}\a_{nr}^s F_{s,k}^{vn}\pem\ .\tag\kn-7$$
We replace the second summation range as follows
\nl
$ (n,s>0)=(s,n>0)+(n>0,s\le 0) -(s>0,n\le 0)\ $
and obtain
$$E_\ep^{(0)}=
\sum_{n>0}\sum_{s}\left(\a_{ns}^v F_{r,k}^{sn}-\a_{nr}^s F_{s,k}^{vn}
\right)\pem
+
\bigg(\sum\Sb n>0\\ s\le 0 \endSb-\sum\Sb s>0 \\ s\le 0 \endSb\bigg)
\a_{nr}^s F_{s,k}^{vn}\pem\ .\tag\kn-8$$
After summation over $s$ in the first sum and using the
``delta distribution''
we see that it will vanish.
Using the integral representation of
$F_{s,k}^{vn}$ (Lemma~\kn.1) and of $\a_{nr}^s$ (\knset-14)
the second part can be rewritten as
$$
\cintd A_r(Q')e_k(Q)\w^v(Q)\bigg(\sum\Sb n>0 \\ s\le 0
 \endSb-\sum\Sb s>0\\
n\le 0 \endSb\bigg)A_n(Q')A_s(Q)\w^s(Q')\w^n(Q)\pem\ .\tag\kn-9$$
Now we use
\proclaim{Lemma \kn.3}
(Bonora et al. \cite\BRRW) For every $N$ we have
$$
\left(\sum_{n>N}\sum_{s\le N}-\sum_{s>N}\sum_{n\le N}\right)
A_n(Q')A_s(Q)\w^s(Q')\w^n(Q)
=d'\Delta(Q',Q) \  .\tag\kn-10
$$
Here $d'$ means differentiation with respect to the variable
$Q'$.
\endproclaim
For completeness we will supply a proof of it below.
Applying Lemma~\kn.3 to our situation we obtain
$$\multline
E_0^{(0)}=
\cintd A_r(Q')e_k(Q)\w^v(Q)d'\Delta(Q',Q)\\
=-\cintd  d'A_r(Q')e_k(Q)\w^v(Q)\Delta(Q',Q)
=-\cintt  dA_r(Q)e_k(Q)\w^v(Q)
=-K_{r,k}^v\ .
\endmultline
$$
If we collect all  the non-vanishing parts
we just obtain the claim of  Proposition~\knsu.1.
\subheading{Proof of Lemma \kn.3}
First we proof the  following relation
\proclaim{Lemma \kn.4}
$$ \gamma_{rk}=(\sum_{n>0}\sum_{s\le 0}- \sum_{s>0}\sum_{n\le 0})
\a_{rn}^s\a_{ks}^n\ .\tag\kn-11$$
\endproclaim
\demo{Proof}
It is the idea of Bonara and collaborators to use representations
via semi-infinite forms. Here we take forms of weight $0$.
Take $\ \Phi=A_1\wedge A_2\wedge \ldots\ $  the vacuum vector
of weight 0 and level 1.  The element $A_i\in \A$
operates with Leibnitz rule as
$$A_i.\Phi=(A_i\cdot A_1)\wedge A_2\wedge\ldots\ +\
A_1\wedge (A_i\cdot A_2)\wedge A_3 \ldots\ +\ \ldots\tag\kn-12$$
As long as $|i|$ is big enough this makes perfect sense.
For some critical strip of indices (e.g\. in particular
for  $A_0=1$) the action has to be modified
\footnote{In physicists' language, it has to be regularized.}
and we obtain only a representation of  a centrally extended algebra
$\Ah'$ where the defining cocycle is  a local one
\cite\KNFA.
Because $\A$ is an abelian Lie algebra two different cocycles
are never cohomologous.  Similar to the case of the
vector field algebra it can be proven that there is  only
one (up to multiplication by a scalar $d$) local cocycle for the
algebra $\A$ (\cite{\KNFA},\ \cite{\KPC}).
Hence
$\ [\widehat{A}_r,\widehat{A}_k]\ldot \Phi=d
\gamma(A_r,A_k)\Phi=d\gamma_{rk}\Phi$.
If $r$ and $k$ are outside of the critical strip of the indices
the action  of $\widehat{A}_r$, resp\. $\widehat{A}_k$ coincides
with the action of the corresponding element $A$ given by (\kn-12).
Inside the critical strip at least
$\  [\widehat{A}_r,\widehat{A}_k]\ldot \Phi\ $ can be
calculated as  follows
(e.g\. see \cite{\SCHLTH, p.137}) for details).
One has to
take only in account  the ways the element $A_s$ inside $\Phi$ will
reproduce itself. First $A_k\cdot A_s=\sum_n \a_{ks}^n A_n$.
This term will only occur if $s\ge 1$ and only the terms with $n<1$
will not be annihilated by neighbouring elements.
To bring it back to $A_s$ by operation of $A_r$ we obtain
$\ \sum_{s>0}\sum_{n\le 0}\a_{rn}^s\a_{ks}^n$.
 Applying the same to $\ -A_k\cdot A_r$,
and changing the variables we obtain
$$[\widehat{A}_r,\widehat{A}_k]\ldot \Phi=
-\big(\sum_{n>0}\sum_{s\le 0}- \sum_{s>0}\sum_{n\le 0}\big)
\a_{rn}^s\a_{ks}^n\Phi=d\gamma_{rk}\Phi\ .\tag\kn-13$$
To determine the constant we calculate this expression for
$r=i$ and  $k=-i$ if $i\gg 0$.  Note that $A_i\Phi=0$, hence
$\ [A_i,A_{-i}]\Phi=A_i(A_{-i}\Phi)$ and that
$A_k\cdot A_s=A_{k+s} + A_j$-terms with indices of $j$ bigger as $k+s$.
Hence we pick up for every $s$ the factor 1  as long as $s-i\le 0$.
Now $s\ge 1$ and we obtain
$\ [A_i,A_{-i}]\Phi=i\cdot\Phi=d\cdot\gamma(A_i,A_{-i})\Phi\ $.
But by calculating residues we get $\ \gamma(A_i,A_{-i})=-i\ $ and
hence the claim.
\qed
\enddemo
\noindent
Note that Bonora et al. also indicate a proof by
direct calculation.

Now
$$dA_k(Q)=\sum_r\beta_{kr}\w^r(Q),\qquad
\beta_{kr}=\cintt dA_k(Q) A_r(Q)=\gamma_{rk}\ ,\tag\kn-14$$
and hence using Lemma~\kn.4
$$ d'\Delta(Q',Q)=\sum_k d'A_k(Q')\w^k(Q)=
\sum_{k,r}
\big(\sum\Sb n>0 \\ s\le 0 \endSb -\sum\Sb s>0 \\ n\le 0 \endSb\big)
\a_{rn}^s\a_{ks}^n\w^k(Q)\w^r(Q')\ . \tag\kn-15$$
Now using that $A_i\w^j=\sum_r\a_{ir}^j\w^r$
we obtain the result for $N=0$.
To get it for general $N$, we just compare the summation ranges with
the $N=0$ range. We obtain
as difference
 $\
 (\sum_{s=1}^N\sum_n-\sum_{n=1}^N\sum_s)
\a_{rn}^s\a_{ks}^n\w^k(Q)\w^r(Q')\ $.
Now each of the partial sums
has only finitely many terms. Hence we can do the summation separately.
But doing this for the summation over $n$ in the  first sum
(after writing the coefficients as integrals and using the
``delta distribution'')
and over $s$ in the second sum we obtain the same value and they will
cancel.
\subheading{Proof of Proposition \knsu.3}
Again we write for $\ep\ne 0$.
$$[L_k(\ep),L_l]=\frac 12\sum_{n,m}\sum_i l_k^{nm}[\nord{u_i(n)u^i(m)}
\, ,\,L_l]
\pem\ .\tag\kn-16$$
As explained at the beginning of this section we can ignore
the normal ordering inside the  above commutators and rewrite
(\kn-16) as
$$
\frac 12\sum_{n,m}\sum_i l_k^{nm}\left(u_i(n)[u^i(m),L_l]
+[u_i(n),L_l]u^i(m)\right)\pem\ .
\tag\kn-17$$
We use Prop.~\knsu.1 to evaluate the commutators
and obtain
$$\frac 12(\ce+\k)\sum_{n,m,v}\sum_i l_k^{nm} \left(
K_{m,l}^v u_i(n)u^i(v)+ K_{n,l}^vu_i(v)u^i(m)\right)\pem\ .\tag\kn-18$$
To make things separately well-defined for $\ep=0$ we have to
rewrite everything again in normal order and pick up commutators for
indices which  are not in normal order.
After evaluating this commutators using Lemma  \knsu.1(3) we can
write the result as  sum $\ A_\ep+B_\ep+C_\ep+D_\ep$ where
they are defined as
$$\align
A_\ep&=\frac 12(\ce+\k)\sum_{n,m,v}\sum_il_k^{nm} K_{m,l}^v
\nord{u_i(n)u^i(v)}\pem\ ,\\
B_\ep&=\frac 12(\ce+\k)\sum_{n,m,v}\sum_il_k^{nm} K_{n,l}^v
\nord{u_i(v)u^i(m)}\pem\ ,\\
C_\ep&=-\frac 12\ce(\ce+\k)\dim\g\sum_{v,m}\sum_{n>v}l_k^{nm} K_{m,l}^v
\gamma_{nv}\pem\ ,\\
D_\ep&=-\frac 12\ce(\ce+\k)\dim\g\sum_{n,m}\sum_{v>m}l_k^{nm} K_{n,l}^v
\gamma_{vm}\pem\ .
\endalign
$$
By considering the range where the coefficients $\ l_k^{nm}\ $ and
$\ K_{m,l}^v\ $ could  be nonzero and taking the normal ordering
into account we see that $A_0$ and $B_0$ are again well-defined
operators and we can ignore the $\pem$ factor.
Renaming  the variables in $B_0$ in the way ($v\to n\to m\to v$) we
obtain
$$A_0+B_0=
\frac 12(\ce+\k)\sum_{n,m,v}\sum_i (l_k^{nm} K_{m,l}^v+
l_k^{mv} K_{m,l}^n)
\nord{u_i(n)u^i(v)}\ .\tag\kn-19$$
The structure constants of the vector field algebra $\L$ can be
calculated as
$$C_{kl}^s=\cintt ([e_k,e_l])\cdot \Omega^s\ .\tag\kn-20$$
\proclaim{Lemma \kn.5}
$$
\sum_m( l_k^{nm} K_{m,l}^v+
l_k^{mv} K_{m,l}^n)=
-\cintt [e_k,e_l]\cdot\w^n\w^v
=
-\sum_s C_{kl}^s l_s^{nv}
\ .\tag\kn-21$$
\endproclaim
\demo{Proof}
We can write the right hand side as
$$-\sum_s\cintd [e_k,e_l](Q)\cdot \Omega^s(Q)\w^n(Q')\w^v(Q')e_s(Q')\ .$$
After summation over $s$ we obtain the ``delta distribution'' for
the pair $(-1,2)$, integrate over $Q'$ and obtain
$\ -\cintt [e_k,e_l](Q)\w^n(Q)\w^v(Q)\ $,
the expression in the middle.
On the left hand side we obtain for the  first sum
$$\sum_m l_k^{nm}K_{m,l}^v=\sum_m\cintd
\w^n(Q)\w^m(Q)e_k(Q)d'A_m(Q')e_l(Q')\w^v(Q')\ .$$
Now applying $\sum_m d'A_m(Q')\w^m(Q)=d'\Delta(Q',Q)$
we obtain after integration over $Q'$
$\  -\cintt \w^n(Q)e_k(Q)d(e_l(Q)\w^v(Q))\ $.
For the second sum we obtain
\nl
$\ -\cintt \w^v(Q)e_k(Q)d(e_l(Q)\w^n(Q))
=\cintt d(\w^v(Q)e_k(Q))e_l(Q)\w^n(Q)\ $.
Together
$$-\cintt \bigg(\w^n(Q)e_k(Q)d\big(e_l(Q)\w^v(Q)\big)
-\w^n(Q)e_l(Q)d\big(e_k(Q)\w^v(Q)\big)\bigg)
\ .$$
If we represent each form by its local representing function
we obtain for the integrand
$$\w^n(z)\w^v(z)\big(e_k(z)\pfz {e_l}(z)-e_l(z)\pfz {e_k}(z)\big)=
\w^n(z)\w^v(z)[e_k,e_l](z)\ .$$
Hence the  claim.\qed
\enddemo
Now
$$A_0+B_0=-\frac 12 (\ce+\k)\sum_s\sum_{n,v}\sum_i C_{kl}^s l_s^{nv}
\nord{u_i(n)u^i(v)}=-(\ce+\k)\sum_s C_{kl}^s L_s\ .$$
It remains to study
$\ \alpha(k,l):=\lim_{\ep\to 0}(C_\ep+D_\ep)$.
Because $L_k$ and $L_l$ are well-defined operators inside $gl(V)$
the scalar  $\a(k,l)$ is  well-defined. Indeed,
using the Jacobi identity inside $gl(V)$ and  the fact
that the $C_{kl}^s$ fulfil
also the Jacobi identity (as they are the structure constants of
$\L$) we obtain that it defines a 2-cocycle.
By studying the order of the forms involved in
defining $l_k^{nm}$, $K_{m,l}^v$ and $\gamma_{nv}$ at the
points $P_+$ and $P_-$ and calculating residues to evaluate
the integral we see that for generic $n,m,v$ values
$$
\gather
 K_{m,l}^v\ne 0\implies -3g\le -v+l+m\le 0,\
\\
l_k^{nm}\ne 0\implies -g\le k-(n+m)\le 0,\quad
 \gamma_{nv}\ne 0\implies -2g\le n+v\le 0\ .
\endgather $$
Adding them up we obtain
$\ \a(k,l)\ne 0\implies -6g\le k+l\le 0$.
For $g\ne 0$ and $n=-g$ and $v=-g-1$ (or vice versa) $\ga_{-g,-g-1}$ will
be non-zero. But in this case the other coefficients will again
yield the same bound.
In particular the cocycle is local.
It can be given as
$$\a(k,l)=-\frac 12 \ce(\ce+\k)\dim \g
\sum_n\bigg(\sum_{m,s}\sum_{v<n}l_k^{nm} l_l^{vs}\gamma_{sm}\gamma_{nv}
+ \sum_{m,s}\sum_{v>m}l_k^{nm} l_l^{vs}\gamma_{sn}\gamma_{vm}\bigg)\ .$$

In the rest of this section we want to find a nicer representation of
the cocycle.
We do not want to take the overall factor $(-1/2) \ce(\ce+\k)\dim \g$
 through
all the calculation. So we just  ignore it in the calculation and
keep it in mind. Let $N$ be a fixed integer.
We start again from $C_\ep$ and $D_\ep$.
We define
$\
E_\ep:=\sum_{v,m}\sum_{n>N}l_{k}^{nm} K_{m,l}^v\gamma_{nv}\pem,\ $
and $\
F_\ep:=\sum_{n,m}\sum_{v>N}l_{k}^{nm} K_{n,l}^v\gamma_{vm}\pem\ $.
In $C_\ep-E_\ep$ and $D_\ep-F_\ep$ only finitely many terms will
occur. Hence we can  put $\ep=0$ and forget about $\pem$.
If we rename in $D_\ep-F_\ep$ the variables $(v\to n\to  m\to v$)
we obtain the following expression
$$\align
\lim_{\ep\to  0}((C_\ep-E_\ep)+(D_\ep-F_\ep))&=
\sum_v\sum_{n=v+1}^N\left(\sum_m l_k^{nm} K_{m,l}^v +
l_k^{mv} K_{ml}^n\right)  \gamma_{nv}
\\
&=-\sum_v\sum_{n=v+1}^N \sum_s C_{kl}^s l_s^{nv}
\gamma_{nv}
\  ,\tag\kn-22
\endalign$$
by applying Lemma~\kn.5.
We  obtain that contrary to the case considered in the proof of
 Prop.~\knsu.1
this will not vanish.
In particular this is the reason that the value of the cocycle will
depend on the normal ordering chosen.

Now we are examining $E_\ep+F_\ep$. We decompose the
summation range for $F_\ep$ as $\ (n,v>N)=(v,n>N)+(v>N,n\le N)
-(n>N,v\le N)\ $
and call the first sum $F_\ep^{(1)}$ and the two others
together  $F_\ep^{(2)}$.
First we obtain
by using
$\gamma_{vm}=-\gamma_{mv}$ and $K_{v,k}^n=\sum_m l_{k}^{nm}
\gamma_{mv}$
$$E_\ep+F_\ep^{(1)}=
\sum_v\sum_{n>N}\bigg(\big(\sum_m l_k^{nm}K_{m,l}^v\big)\gamma_{nv}-
K_{n,l}^v K_{v,k}^n\bigg)\pem\ .$$
Because  $F_\ep^{(2)}$ will be well-defined  for $\ep=0$ (see below)
 we can
ignore here $\pem$ as  long as we not break it into the two partial sums.
Using the integral representation of the coefficients
and performing the $m-$summation we get a ``delta distribution''
(as we did it several times) and obtain
$$\align
\bigg(\sum_m l_k^{nm}K_{m,l}^v\bigg)\gamma_{nv}&=
-\cintd w^v(Q)d'A_n(Q')e_l(Q)d\big(e_k(Q)\w^n(Q)\big)A_v(Q')\ ,\\
K_{n,l}^v K_{v,k}^n&=
\cintd   w^v(Q')d'A_n(Q')e_l(Q')d(e_k(Q)\w^n(Q))A_v(Q)\ .
\endalign
$$
Summation over $v$ and integration over $Q'$ shows
that $E_\ep+F_\ep^{(1)}=0$.

It remains to look at
$$F_\ep ^{(2)}=
\bigg(\sum_{n>N}\sum_{v\le N}-\sum_{n\le  N}\sum_{v> N}\bigg)
K_{v,k}^n K_{n,l}^v\,\pem\ .\tag\kn-23$$
By checking the structure constants we see that every partial sum
consists only of finitely many terms. Hence we can set again $\ep=0$.
For $N=0$
we collect the surviving terms (\kn-22) and (\kn-23) and
obtain the result of the structure of the cocycle (\knsu-11) as
claimed in Prop.~\knsu.3.
If we choose a different normal ordering we would just  obtain
a different summation prescription in  (\kn-22). The expression (\kn-23)
would be the same.
It remains to show the form (\knsu-12) of the cocycle for the
special basis elements. We look first at
$\widehat{\chi}_{k,-k}$. Now
$K_{v,k}^n\ne 0$ implies $k\le n-v$ and
$K_{n,-k}^v\ne 0$ implies $-k\le v-n$, hence $v=n-k$. And as value
we obtain in this case $(n-k)$ resp\. $n$.
Assume $k\ge 0$ then only the first sum will survive.
It will be the finite sum
$$
\sum_{n=1}^k (n-k)n=-\frac 16(k^3-k)\ .$$
In the expression for $\psi_{k,-k}$ we obtain that
$\ C_{k,-k}^s l_s^{nv}\gamma_{nv}\ne 0$ only if
$s\ge 0$, $s\le n+v$ and $n+v\le 0$. This implies
that we obtain non-zero terms only for $s=0$ and $n=-v$.
But if we look at the summation range $\ \sum_v\sum_{n=0}^{v+1}\ $
we see that $n=-v$ imposes that either $v=0$ (with vanishing coefficient)
or that the summation range is empty.
Hence $\psi_{k,-k}=0$.
Altogether we get the form (\knsu-12).
\vskip 1cm
For further reference let us note
\proclaim{Corollary~\knpr.1}
Let $E$ be an operator of $gl(V)$, $V$ the fixed admissible representation
chosen in Section \knsu,  such that there exists a basis element $e_l$ of
$\L$ with
$$ [E,x(r)]=-(\ce+\k)\,x(\nabla_{e_l}A_r)\ ,\tag\kn-24
$$ for every $x\in\g$ and every $r$,
then for every $k$ we have
$$[L_k,E]=[L_k,L_l]\, $$
where the latter can be evaluated by (\knsu-10).
\endproclaim
\demo{Proof}
In the proof of Prop.~\knsu.3 we used only the relation
$\ [L_l,x(r)]=-(\ce+\k)\,x\otimes(\nabla_{e_l}A_r)\ $ of Prop.~\knsu.2.
Hence if $E$ is an operator obeying the same rules we will
obtain
$\ [L_k(\ep),E]=[L_k(\ep),L_l]\ $ and hence also in the limit $\ep\to 0$:
$\ [L_k,E]=[L_k,L_l]$.
\qed
\enddemo
%
\vskip 1cm
\head
5. The weight of the Sugawara representation
\endhead
\def\kn{5}
\vskip 0.2cm
Assume that we have a $\Gh$-module with a certain highest weight
such that $\Lh$ acts via the  corresponding Sugawara representation
on it as introduced in Section 5.
The goal of this
section is to express the weight of the Sugawara representation in
terms of the weight of the $\Gh$-module under consideration. The
following
result will generalize the result of \cite\SHEINNS , which was
obtained for the \KN\ algebras of Heisenberg type.

The theory of highest weight modules for \KN\ algebras of affine type
is developed in \cite\SHEINW\ (see also \cite\SHEINNS ). Let us
outline the class of highest weight modules we are considering  here.

Let $\g$ be again a fixed finite-dimensional simple Lie algebra.
We denote by $\frak h$ a fixed Cartan subalgebra, by $\frak n_+$ (resp\.
$\frak n_-$) its upper (resp\. lower) nilpotent algebra, and by
$\frak h^*$ the dual of the Cartan subalgebra.
As usual $U(B)$ denotes the universal enveloping algebra of the
Lie algebra $B$.
Note that $\g$ is embedded via
$\ x\mapsto x\otimes A_0\ $ into $\G_0$ and $\Gh_0$. Hence
$\frak h,\; \frak n_+,\;\frak n_-$ can  be considered as subalgebras
of $\Gh$.

Let $V$ be an almost-graded ${\Gh}$-module.
Recall from Section 2 that $\Gh$ is an almost-graded Lie algebra.
Let us recall the definition of an {\it almost-graded module}.
We have given a $\Z$-degree deg on $V$. The elements
$x$ with $\ \deg x=n\ $ are
called homogeneous elements of degree $n$. Denote by $V_n$ the
subspace generated by all homogeneous elements of degree $n$.
It is assumed to be finite dimensional
\footnote{
For certain applications the finite dimensionality of the
homogeneous subspaces might be dropped. With respect of
this more general conditions our modules are what is called
``{\it quasi-finite modules}''.
}
 for every $n$ and $V$ is
the vector space direct sum $\ V=\sum_{n\in\Z} V_n$.
Concerning the action of $\Gh$ on $V$, there are
constants $K,L$ such that for all $n$ and $m$ we have
$$\Gh_m\ldot V_n\ \subseteq\ \sum_{h=n+m-K}^{n+m+L} V_h\ .$$

Because we are heading for the Sugawara representation we will assume
that a certain normal ordering $\nord{...}$ has been fixed.
 We will suppose that $V$ together with the normal ordering
satisfy the following conditions:
 \roster
  \item"{(1)}"  There exists an element $v\in V$ such that
    $$\Gh _+v={\frak n}_+v=0,\quad\text{and}\quad
       V=U(\Gh )v.                                    \tag 5-1
      $$
  \item"{(2)}"
      There exists $\ce\in\C$ and
    $\chi =(\chi_0,\chi_{-1},\ldots ,\chi_{-g})\in(\frak h^*)^{g+1}$,
    such  that for all $\ m,n\in\{0,-1,\ldots ,-g\}$ and
    all $\ h,h'\in\frak h\ $ we have
     $$ h(m)v=\chi _m(h) v +\Tld,\quad
     \text{and}\quad t.v=\ce v,          \tag 5-2
     $$
    and
     $$ \nord{h(m)h^\prime (n)}v=
        \chi _m(h)\chi _n(h^\prime )v +\Tld \ .           \tag 5-3
     $$
     Here and in the following $\Tld$ denotes terms of lower degree.
     Recall that $\ t\ $ is the basis element of the center.
     Such a tuple $\chi$ is called the {\it weight},
       $\ce$ is called the {\it central charge},
     and $v=v_{\chi,\ce}$  is called the {\it highest weight vector}
    of the module $V$.
  \item"{(3)}"  For any positive finite root $\a$
     and corresponding root vector $\ea_\a$  we have for all
    $m=0,-1,\ldots ,-g$,
     $$ \ea_\a (m)v=0\cdot v+\Tld\  .                            \tag 5-4
     $$
  \item"{(4)}"
     The almost-grading of the $\Gh$-module $V$ is compatible
     with the almost-grading of the associative algebra of functions
     $\A$ in the following sense.
     For $n,m\le -g$ and every $v_m\in V$ with $\deg v_m=m$ we have
     that in $u(n)v_m$ elements of degree $k$ could only occur if
     in $A_nA_m$ elements of degree $k$ are occuring.
     In formulas, if $\{v_k^1,v_k^2,\ldots,v_k^r\}$ are   basises of the
     spaces $V_k$  then
   $$
    \multline\qquad\qquad\{k\mid
   u(n)v_m=\sum_{k,s}d^{k,s}_{nm}(u)v_k^s ,\ \exists s:d^{k,s}_{nm}(u)\ne
    0\}
    \\
    \subseteq
    \{k\mid  A_nA_m
    =\sum\limits_k\a^k_{nm}A_k,\ \exists \a^k_{nm}\ne 0\}\ .
    \endmultline $$
 \item"{(5)}"
    In     addition we suppose that ${\deg v}<-g$
    and that all elements of degree $\ge\deg v$ are
    multiples of $v$.
 \endroster
Conditions (4) and (5) will take care that for $n<-g$
the vector $\ x(n)v\ $ in its decomposition  will
contain only elements of degree less than the degree of $v$.
\bigskip
Now we are coming to the Lie algebra $\Lh$. We call a tuple
$\ \l=(\l_0,\l_{-1},\ldots,\l_{-3g})\in\C^{3g+1}$ a {\it weight} of the
Lie algebra $\Lh$.
To each weight $\lambda$ we assign the quadratic differential
    $$\Lambda =\sum\limits_{k=-3g}^{0} {\lambda}_k{\Omega}^k, \tag 5-5
    $$
where the $\Omega^k$ are  the basis elements
introduced in Section \knset, see (\knset-11).
This quadratic differential
$\Lambda$ will be called a {\it weight} of the Lie algebra $\Lh$ as well.
Note that also for the $\Gh$-modules (resp. for the  $\Ghe$-modules
which we are considering in the next section)
 the  weights  can be
identified with the abelian differentials of first order
\cite{\SHEINW,\SHEINNS}.

Let  $V$ be an $\Lh$-module.
In analogy with (1)-(5) we  impose the following conditions on the
module $V$.
 \roster
  \item"{(6)}"  There is an element $v\in V$ such that
    $$\Lh _+v= 0,\quad\text{and}\quad V=U(\Lh )v\ .\tag 5-6
    $$
  \item"{(7)}"  There exists a weight $\lambda$  such  that
     $$ E_mv=\lambda _m v +\Tld ,                            \tag 5-7
     $$
    for each $m=0,-1,\ldots ,-3g$.
    Here we denote by $E_m$ the element $(e_m,0)$ which is the
    lift of the vector field $e_m\in\L$ to the central extension
    $\Lh$.
 \endroster
\bigskip
Now let us define  more precisely which classes of normal orderings
$\Sigma$ we
consider in this section.  Let ${\Sigma}^{\pm}$  be  a
decomposition of the two-dimensional lattice
\nl $\{\;(m,n)\mid m,n\in\Z\;\}$ into
two disjoint parts
such that the ${\Sigma}^{\pm}$ differ from
$${\Sigma}_0^{\pm}=\{(m,n)\mid m\le n\,\,(\text{resp\.\;}m > n),\
  m,n\in\Z\}$$
 only for $-g\le m,n\le 0$.
 According to \cite\KNFA\  each decomposition of
this kind defines a normal ordering by setting
     $$\nord{x(m)y(n)} \quad =
           \ \cases x(m)y(n)&,(m,n)\in{\Sigma}^+\\
                    y(n)x(m)&,(m,n)\in{\Sigma}^-\ .
             \endcases                                     \tag 5-8
     $$
for each $x,y\in\fg ,$ and each pair $(m,n).$
\nl
More general normal orderings with the only restriction
that they differ from $\Sigma_0$ at most for finitely many
pairs $(m,n)$ are possible.
There is no essential difficulty involved.
It is only that the equations in the following would look
slightly more complicated.
\remark{Warning}
For abelian $\g$ by the decomposition $\Sigma^\pm$
the order of the elements inside the normal ordering colons does not
play any role, e.g\. $\ \nord{x(n)y(m)}=\nord{y(m)x(n)}$.
For non-abelian $\g$ this is only true for the elements of
the Cartan subalgebra $\frak h$.
\endremark
\noindent
For further reference let us denote
$\Sigma_{cs}^{\pm}:=\Sigma^{\pm}\cap\{(m,n)\mid
 -g\le m,n\le 0,\ m,n\in\Z\}$.
\bigskip
In the rest of this section we would like to connect $\chi$ and
$\lambda$ assuming that $V$ is a $\Gh$-module of  highest weight
$\chi$ and $\lambda$ is the weight of the corresponding Sugawara
representation.
Assume that $V$  obeys the Conditions
(1)--(5). Let $\chi$ be its weight, $\ce$ its central charge and
$v$ its highest weight vector.
By (\knwei-1) we see that it is an admissible module in the sense of
Section 3. Hence we can apply the Sugawara construction. Assume
$(\ce+\ka)\ne 0$  (which is for positive $\ce$ always the case),
then we obtain by mapping $E_k$ to
$\ L_k^*=\dfrac {-1}{\ce+\ka}L_k$, the rescaled modes of the
Sugawara operator,
and mapping the central element of $\Lh$ to the identity operator on $V$
a representation of $\Lh$. This is the content of Theorem \knsu.1.
\proclaim{Lemma 5.1}
Independently of the normal ordering chosen
for this representation of $\Lh$ we have
$$L_k^* v=0,\quad\text{for all\ } k >0,\tag\kn-9$$
and
$$L_k^* v=\l_k v+\Tld,\quad\text{for \ }
   k=0,-1,\ldots, -3g,\ \tag\kn-10$$
with certain $\l_k\in\C$.
\nl
The elements $\ L_k^*v\ $ for $k < -3g$ contain in their degree
decomposition only elements $w$ with $\ \deg w <\deg v$.
\nl
In particular the Condition (6) and (7) are fulfilled
for the subspace $\ U(\Lh)v$.
\endproclaim
\demo{Proof}
By definition $L_k^*v$ is up to a scalar given as
$$w:=\sum_{n,m}\sum_{i} l_k^{mn}\nord{u_i(m)u^i(n)}v\ .$$
Note that due to the summation over $i$ (see Lemma \knsu.1(3))
we obtain
$$
\sum_i  l_k^{mn}u_i(m)u^i(n)=
\sum_i  l_k^{mn}u_i(m)u^i(n)-\dim\g\,\ga_{mn}l_k^{mn}\ .$$
But $\ \ga_{mn}l_k^{mn}\ne 0\ $ implies $-3g\le k\le 0$.
Hence for $k$ outside of this range the order of $(m,n)$ is of
no importance.
In particular we can always take the order corresponding to the
normal ordering $\Sigma_0$.
Note also that $l_k^{nm}\ne 0$ implies $-g\le k-(n+m)\le 0$.
Hence if $k>0$ either $n$ or $m>0$ for the  terms
where $l_k^{nm}\ne 0$.
 By this and (\kn-1) we get
$L_k^*v=0$ in this case, hence (\kn-9).
For $k<-3g$ we see that only pairs $(n,m)$ with
$n,m\le 0$ will occur.
Again $l_k^{nm}\ne 0$  implies that
 $(n+m)<-2g$ which could only be the
case if $n$ or $m$ is less than $-g$.
Now (4), (5)
and the almost-grading of the algebra $A$
implies that the resulting vector
$\sum_i u_i(m)u^i(n)v$ will have only elements of degree
$<\deg v$.
Hence the claim for $L_k^*v$ for $k<-3g$. The equation
(\kn-10) is automatic because $v$ is the element of
highest degree.
Note that $\nord{...}$ was allowed to be an
arbitrary normal ordering.
\qed
\enddemo
Using the arguments in the proof above we see that with
respect to the normal ordering (5-8) the $\deg v$-part
of $L_k^*$ for $-3g\le k\le 0$ will only come from
$$
L_k^*\, v=\frac 12\sum_{-g\le m,n\le 0}l_k^{mn}\nord{u_i(m)u^i(n)}v+\Tld\ .
$$
For more general normal orderings an additional term $M(k,\Sigma)v$ will
occur where  $M(k,\Sigma)$ is a scalar depending on the
normal ordering.

By definition $\l=(\l_0,\l_{-1},\ldots,\l_{-3g})$
with the $\l_i$ from Lemma~\kn.1
is the weight of the Sugawara representation.
As usual we denote by $\Lambda$ the corresponding quadratic differential.
\proclaim{Lemma 5.2}
Let $T$ be the Sugawara operator (or energy-momentum tensor) as
introduced by (\knsu-3) then
$$T.v=-(\ce+\ka)\,\Lambda\cdot v+\Tld\ .\tag\kn-11
$$
\endproclaim
Note that $T.v$ is a formal sum. But if we decompose it into its
homogeneous parts  with respect to the degree in $V$ then the
components are well-defined elements
of $V\otimes\Fn 2$ ($\Fn 2$ is the space of quadratic differentials).
The Equation (\kn-11) should be interpreted in this sense.
\demo{Proof of Lemma \kn.2}
$$T.v=\sum_{k\in\Z}\Omega^kL_kv=\sum_{k\le 0}\Omega^kL_kv\ .$$
By isolating the highest degree part using Lemma \kn.1 we see
$$(T.v)_{h.D.}=-(\ce+\ka)\sum_{k=-3g}^0\Omega^kL_k^*v
=-(\ce+\ka)\sum_{k=-3g}^0\Omega^k\l_k v
=-(\ce+\ka)\Lambda\cdot v\ .
\qed$$
\enddemo

The following theorem answers the question which was put in the
beginning of this section. The analogue of this theorem for the less
complicated case of Heisenberg type algebras (which
corresponds to abelian
$\fg$) was proved in \cite\SHEINNS .
\proclaim{Theorem 5.1} Let $\chi$ be the highest weight of the
 $\Gh$-module $V$ fulfilling the condition (1)--(5), $\Sigma$ be a
 normal ordering and $\Lambda $ be the weight of the corresponding
 Sugawara representation of $\Lh$.
 Then the following relations for $\Lambda$ hold:

    $$\Lambda = \frac {-1}{\ce+\ka}\sum\limits_{-g\le m,n\le 0}
         \bigg(\frac1{2}<\chi _m,\chi _n>+
         \sum_{-g\le s\le 0}\a_{mn}^s<\chi _s,{\bar\rho}>-
         \npos\cdot\ce\,\ga_{[mn]} \bigg)\w ^m\w ^n
                                                      \tag 5-12
    $$
 and
    $$\Lambda = \frac {-1}{\ce+\ka}\sum\limits_{-g\le m,n\le 0}
         \bigg(\frac1{2}<\chi _m,\chi _n>+
         \sum_{-g\le s\le 0}\a_{mn}^s<\chi _s,{\bar\rho}>
         \bigg)\w ^m\w ^n+\Ksig,
                                                      \tag 5-13
    $$
with
$$\Ksig:=\frac {\npos\cdot\ce}{\ce+\ka}
\cdot\big(\sum_{(m,n)\in\Sigma_{cs}^+} -
\sum_{(m,n)\in\Sigma_{cs}^-}\big)\ga_{mn}\w^m\w^n\ .$$
Here the $\a_{mn}^s$ are
the structure constants of the associative algebra $\A$
defined by (2-13), $\ga_{[mn]}:=\pm\ga_{mn}$ if $(m,n)\in\Sigma^\pm$
respectively, $\ 2{\bar\rho}=\sum_{\a >0}\a\  $ is
the sum over the positive roots,
and $\npos$ is the number of positive roots, i.e\.
$\npos=1/2(\dim\g-\rank\g)$.
\nl
For a normal ordering $\Sigma$ with
$\Sigma^+_{cs}$ (or $\Sigma^-_{cs}$) $\subseteq
\{(m,n)\mid -g\le m,n\le 0\}\ $ we obtain $K(\Sigma)=0$.
\endproclaim
\demo{Proof} By Lemma 5.2 it is enough to find the term of highest
degree in $T.v$, where $T$ is given by (3-2). Note that by (4),(5)
$\nord{u_i(m)u^i(n)}v$ has a nonzero projection onto ${\Bbb C}\,v$ only if
$-g\le n,m\le 0$. So
    $$ T.v=\bigg({1\over 2}\sum\limits_{-g\le m,n\le 0}\sum_i
      \nord{u_i(m)u^i(n)}{\omega}^m{\omega}^n\bigg)v + \ldots
                                                            \tag 5-14
    $$
For the rest of the proof let $-g\le m,n\le 0$.
We choose as basis in $\g$ the ``canonical basis''. Choose
$h_k,\ k=1,\ldots,\rank \g$ for the Cartan subalgebra and
$h^k,\ k=1,\ldots,\rank \g$ as dual basis with respect to
the invariant form $(..|..)$. Choose for the positive roots
$\a$ basis elements $\ea_\a$ and take as dual $\ea_{-\a}$
the rescaled element of the  corresponding negative root
(e.g\. $(\ea_\a|\ea_{-\a})=1$).
With respect to this basis $\ (h_k,\;\ea_\a,\;\ea_{-\a})$
   $$\sum_i u_i(m)u^i(n)=\sum\limits_{\a >0}
       \bigg(\ea_\a (m)\ea_{-\a}(n)+\ea_{-\a}(m)\ea_\a (n)\bigg)
       +\sum\limits_kh_k(m)h^k(n),
                                                            \tag 5-15
   $$
where $k=1,\ldots , \rank \fg$  and  $\a$ runs over all
positive roots.
Using the structure equation of $\Gh$ we can write
   $$ \multline
     \ea_\a (m)\ea_{-\a}(n)=\ea_{-\a}(n)\ea_\a (m)+[\ea_\a (m),
        \ea_{-\a}(n)]= \\
        =\ea_{-\a}(n)\ea_\a (m)+\sum_s\a_{mn}^sh_\a (s)
        -\gamma_{mn}\, t\ ,
      \endmultline                                        \tag 5-16
   $$
where $h_\a:=[\ea_\a ,\ea_{-\a}]$ and the $\ga_{mn}$ are defined by (3-8).
We replace the summands in (\kn-15) at the first position
by (\kn-16). Now  $(\ea_{-\a}(m)\ea_\a (n))v=0\cdot v+\ldots $, by (\kn-4),
hence there remains only
   $$\sum_i\ u_i(m)u^i(n) v=
      \bigg(\sum\limits_k h_k(m)h^k(n) +
       \sum_{s}\sum\limits_{\a >0}\a_{mn}^sh_\a (s)
       -\npos\cdot\gamma_{mn}t\bigg)v +\ldots .
                                                            \tag 5-17
$$
If $\ (m,n)\in\Sigma_{cs}^+\ $ then this is already
in normal order.
Note that $h_\a (s)v$ will have a $v$-component only if $-g\le s\le 0$.
Now by (5-2), (5-3) we obtain
   $$\bigg(\sum\limits_k\chi_m(h_k)\chi_n(h^k)
        +\sum_{s=-g}^0\sum\limits_{\a >0}\a_{mn}^s\chi_s(h_\a )
        -\npos\cdot\gamma_{mn}\ce\bigg)v+\ldots \ .
                                                            \tag 5-18
   $$
For the pairs $(m,n)\in\Sigma_{cs}^-$ we get
 $$
\aligned\sum_i\nord{u_i(m)u^i(n)}v&=\sum_i u^i(n)u_i(m)v
\\
&=
      \bigg(\sum\limits_k\nord{h_k(m)h^k(n)}+
       \sum_{s}\sum\limits_{\a >0}\a_{nm}^sh_\a (s)
       -\npos\cdot\gamma_{nm}t\bigg)v +\ldots
\\
&= \bigg(\sum\limits_k\chi_m(h_k)\chi_n(h^k)
        +\sum_{s=-g}^0\sum\limits_{\a >0}\a_{mn}^s\chi_s(h_\a )
        +\npos\cdot\gamma_{mn}\ce\bigg)v+\ldots \ ,
\endaligned                               \tag 5-19
$$
where we used $\a_{mn}^s=\a_{nm}^s$ and $\ga_{mn}=-\ga_{nm}$.
\nl
Let us denote the dual bilinear form on $\frak h^*$ by
$<..,..>$.
We use
 $$\ \sum_k\chi_m(h_k)\chi_n(h^k)=<\chi_m,\chi_n>,\qquad
 \chi_s(h_\a)=<\chi_s,\a>, \quad\text{and}\quad
    \sum\limits_{\a >0}\a =:2{\bar\rho},
$$ in (\kn-18) and (\kn-19)
and obtain after summation over $n$ and $m$
   $$\split
      Tv &=\bigg({1\over 2}\sum\limits_{-g\le m,n\le 0}\sum_i
      \nord{u_i(m)u^i(n)}{\omega}^m{\omega}^n\bigg)v + \ldots  \\
      &=  \sum\limits_{-g\le m,n\le 0}\bigg(\frac1{2}<\chi _m,\chi _n>+
         \sum_{-g\le s\le 0}\a_{mn}^s<\chi _s,{\bar\rho}>-
         \npos\cdot\ce\,\ga_{[mn]} \bigg)\w ^m\w ^n v+\Tld  \\
      &= \sum\limits_{-g\le m,n\le 0}\bigg(\frac{1}{2}<\chi_m,\chi_n>+
        \sum_{s=-g}^0\a_{mn}^s<\chi_s,{\bar\rho}>
        \bigg){\omega}^m{\omega}^n v-(\ce+\ka)\Ksig v
        +\ldots .
     \endsplit                                          \tag 5-20
   $$

After dividing this by $-(\ce+\ka)$ we obtain equations (\kn-12)
and (\kn-13).
For the special normal ordering $\Sigma$ in the theorem
we have $\Sigma_{cs}^-=\emptyset$ and $\Sigma_{cs}^+$ is the
full square. For $(m,n)$  with $m\ne n$ also the
pair $(n,m)$ will appear. Because $\ga_{nm}=-\ga_{mn}$ they
will cancel. $\ga_{mm}$ will always be zero.
Hence, $\Ksig=0$.
The same is true if we interchange the role of $\Sigma_{cs}^+$ and
$\Sigma_{cs}^-$. \qed
\enddemo
Note again that a more general normal ordering will yield that
the summation prescription in the definition of $\Ksig$ will be
different.
\example{Example} Let $g=0.$ Then (5-12) specializes to
    $$\Lambda = \frac {-1}{\ce+\ka}(\frac1{2}<\chi _0,\chi _0>+
         <\chi _0,{\bar\rho}>)(\w ^0)^2,
    $$
  where $\w^0=\frac 1z\text{d}z$.
Note that in this case $\Ksig$ will always vanish.
 Up to the factor $\frac {-1}{2(\ce+\ka)}(\omega^0)^2$
this
  expression coincides with the eigenvalue of Casimir operator of second
  order for $\fg .$
\endexample
%
%
\vskip 1cm
\head
6. The Casimir operator and its eigenvalues
\endhead
\def\kn{6}
\vskip 0.2cm
As it is pointed out in \cite{\KACRAI, Lecture 10}\ the Sugawara
construction is closely related to the Casimir operators of
\KM\ algebras. In this section generalizing the approach of
\cite\KACRAI\ we obtain a Casimir operator of the second order and
under certain assumptions its
   eigenvalues for an arbitrary \KN\ algebra of affine type.

In order to do this we first extend our algebra $\Gh$ by
adjoining an arbitrary vector field $e\in\L$. As vector space we
take
$\ \Ghe=\Gh\oplus\C\cdot e\ $ and define
   $$[e,x\otimes A_n]=x\otimes \nabla_{e}A_n,\qquad [\,e,t]=0\ . \tag\kn-1
   $$
\proclaim{Proposition~\kn.1}
\quad  $\Ghe$ is a Lie algebra.
\endproclaim
\demo{Proof}
We have to show the Jacobi identity for commutators in
which the new element $e$ is involved, i.e\. for the
triples $(e,x(n),e)$, here it is clear, and $(e,x(n),y(m))$, here we
shall show it.
Now
$$\align
[e,[x(n),y(m)]]&=[e\,,\,[x,y]\otimes A_nA_m-(x|y)\gamma(A_n,A_m) t]=
[x,y]\otimes \nabla_e (A_nA_m)
\\
&=
[x,y]\otimes \big((\nabla_e A_n)A_m)+A_n(\nabla_eA_m)\big),
\\
[x(n),[y(m),e]]&=-[x(n),y\otimes \nabla_eA_m]=
-[x,y]\otimes A_n\nabla_eA_m+(x|y)\gamma(A_n,\nabla_eA_m)t,
\\
[y(m),[e,x(n)]]&=[y(m),x\otimes \nabla_eA_n]=
-[x,y]\otimes A_m\nabla_eA_n-(y|x)\gamma(A_m,\nabla_eA_n)t\ .
\endalign
$$
We have to show that the central terms will cancel:
$$\multline
\gamma(A_n,\nabla_eA_m)=\cintt A_nd(\nabla_e(A_m))=
\cintt A_n\nabla_e(dA_m)
\\
=
\cintt \nabla_e(A_ndA_m)-\cintt (\nabla_e A_n)dA_m
=\cintt A_md(\nabla_e(A_n))=\gamma(A_m,\nabla_eA_n)\ .
\endmultline$$
Here we used  that the Lie derivative commutes with the exterior
differentiation, that it is a derivation, and that
the Lie derivative of a meromorphic form has no residue
(see \cite{\SCHLTH, p.102}). \qed
\enddemo
\noindent
Note that in the situation of zero genus \cite\KIL\ one takes
$e=z\fpz=d$ and obtains for (\kn-1)
\nl $\ [e,x(n)]=n\cdot x(n)$.

Let the vector field $e$ be fixed and decomposed as $e=\sum_k\eps^k e_k$
(a finite sum).  Let $V$ be a highest weight representation of $\Gh$
which extends to a highest weight representation of $\Ghe .$
In what follows we do not distinguish
between $e$ as  element of the Lie algebra $\Ghe$ and the
corresponding
operator on the representation space.
The weight $\lambda_e$ corresponding to $e$ is given by
$\ e\ldot v=\lambda_e v+\Tld$.
\proclaim{Theorem~\kn.1}
  Let $\ e=\sum_k\eps^ke_k\ $ be the  vector field used to define
  the algebra $\Ghe$. Let the finite dimensional algebra
  $\g$ be either simple or abelian.
  Let $V$ be a highest weight representation of $\Ghe$ and
  $\ L=\sum_k\eps^kL_k\ $ be the corresponding linear combination of
  the Sugawara operators $L_k$ introduced in Section 3, then
  $$\Omega:=2L+2(\ce+\k)e\tag\kn-2$$
  is a Casimir operator of $\ \Ghe\ $ for the representation  $V$.
  The number $\ce$ is the central charge (i.e\. $t. v=\ce\cdot v$) and
  $\k$ is the dual Coxeter number
  (resp\. $0$ in the abelian case)
  introduced in Lemma~\knsu.1.
\endproclaim
\demo{Proof}
We have to show that $\ \Omega\ $ commutes with all other operators of
$\Ghe$.
\nl
(1)
Note that we have $\ [e,x(n)]=x(\nabla_eA_n)$. On the other hand
by Prop~\knsu.2(1)
$$
[L,x(n)]=-(\ce+\k)\sum_k\eps^k[L_k,x(n)]=
-(\ce+\k)\sum_k\eps^kx(\nabla_{e_k}A_n)=
-(\ce+\k)x(\nabla_{e}A_n)\ .$$
Hence $\ [\Omega,x(n)]=0.$
\nl
(2) It remains to show $[\Omega,e]=0$.
First we want to calculate $[e,L_k]$.
Due to the linearity of the Lie derivative and (\kn-1) the
operator $\ E:=-(\ce+\ka)e\ $ obeys the relation
$$\ [E,x(n)]=-(\ce+\ka)x(\nabla_eA_n)=[L,x(n)]\ ,$$
by Prop.~\knsu.2(1). Hence from Corollary \knpr.1 we get
$\ [L_k,E]=[L_k,L]\ $.
Now
$$-(\ce+\ka)[L,e]=-(\ce+\ka)\sum_k\eps^k[L_k,e]=\sum_k\eps^k[L_k,L]=[L,L]
=0\ .$$
Because $(\ce+\ka)\ne 0$ this implies $[L,e]=0$
and further  $[\Omega,e]=0$.\qed
\enddemo
In the remainder of the section let us discuss the question of
eigenvalues of the Casimir operators. In zero genus case a Casimir operator
acts as  multiplication by a scalar on the  highest weight
module. The usual argument is that the highest weight vector $v$ is an
eigenvector of the  Casimir operator: $\Omega v=\l v$
 \cite\KACRAI\ . Each other
element $w$ of this module can be represented in the form $w=uv$
where $u$ belongs to the universal enveloping algebra (cf.(5-1)).
As $\Omega$ commutes with $u$ one has $\Omega w=\Omega uv=
u\Omega v=\l uv=\l w.$
This argument does not work in the almost-graded
 situation. In general
the  question whether $\Omega$
has even at least one eigenvector
has no obvious answer.
Indeed in the case that the highest weight
vector $v$ is an eigenvector of the Casimir operator $\Omega$ one can
conclude as above that $\Omega $ operates as
 a scalar. This scalar can be found
in the following way. Recall the definition of the  complex number $\l_e$ by
the relation $\ e.v=\l_ev+\ldots\ $.
\proclaim{Proposition 6.1} Let $\eps ^k=0$ for $k<-g$ or $k >0$.
  Suppose that the Casimir
  operator $\Omega$ as defined in (\kn-2) acts in a
  highest weight representation of  weight
  $\chi =(\chi _0,\chi_{-1},\ldots ,\chi_{-g})$ as  multiplication
  by a constant $\lambda_\Omega .$ Then  this constant equals
  $$\lambda_\Omega =\sum\limits_{k,m,n}\eps ^kl_k^{mn}
    \big(<\chi _m,\chi _n>+2\sum\limits_s\a_{mn}^s<\chi _s,{\bar\rho}>
    -2\npos\cdot\ce\,\ga_{[mn]}\big)
    +2(c+\k)\l _e\ ,                                      \tag 6-3
  $$
   (where all summations are over the range $\ -g,\ldots,-1,0$).
The term with $\ga_{[mn]}$ will vanish if a normal ordering
$\Sigma$ is chosen such that $\Ksig=0$.
\endproclaim
\demo{Proof} It follows from the Theorem 5.1 that the following relation
for the highest weight of the Sugawara representation holds:
    $$-(\ce+\ka)\cdot 2\Lambda = \sum\limits_{m,n}\bigg(<\chi _m,\chi _n>+
               2\sum\limits_s\a_{mn}^s<\chi_s,{\bar\rho}>
              -2\npos\cdot\ce\,\ga_{[mn]}\bigg)
               \w ^m\w ^n.
                                                          \tag 6-4
    $$
On the left-hand side we use the definition (\knwei-5)
$\ \Lambda=\sum_{k=-3g}^0\l_k\Omega^k\ $ and on the right hand
side  $\ \w ^m\w ^n=\sum l^{mn}_k\Omega^k$ (see (3-4)).
If we compare both sides we obtain
    $$-(\ce+\ka)2 \l _k=\sum\limits_{m,n}l_k^{mn}\bigg(<\chi _m,\chi _n>+
             2\sum\limits_s\a_{mn}^s<\chi _s,\bar\rho >
             -2\npos\cdot\ce\,\ga_{[mn]}\bigg)\ .
                                                           \tag 6-5
    $$
{}From the definition of $\ \Omega=2L+2(\ce+\ka)e\ $ with
$\ L=\sum_k\eps^kL_k=-(\ce+\ka)\sum_k\eps^k L_k^*$ we obtain
$\Omega v=\lambda_\Omega v$, where $\lambda_\Omega$
is given as the expression (6-3).
\qed
\enddemo
\noindent Again, for more general normal orderings $\Sigma$ an explicit
correction term can be given.
\example{Example} Let $g=0$.
Then the only admissible values of $k,m,n$ are $k=m=n=0,$ and
we have $l^{00}_0=\a_{00}^0=1$.
 Set $e=z\fpz=d\ $ (i.e\. $\eps^0=1$) then the result of
Theorem 6.2 specializes to
  $$ \lambda_\Omega = <\chi_0+2{\bar\rho},\chi_0>+2(c+\k)\l_e\ .
                       \tag\kn-6
  $$
This is the as expected the result, as it can be found (for example) in
\cite{\KACRAI, Prop.10.2}\ $\ \lambda_\Omega = <{\tilde\chi}+2\rho
,{\tilde\chi}>,$ where ${\tilde\chi}$ is the weight of the module
 $V$ (including
the central charge and the vector field d)
 and $\rho$ is the sum of fundamental weights for the
affine Kac-Moody algebra  $\Ghe $.
\endexample
%
%
%
\vskip 1cm
\head
Appendix: Sugawara construction for the multi-point
situation
\endhead
\def\kn{\text{A}}
\vskip 0.5cm
The Sugawara construction above can be generalized to the
situation where one allows poles  at more than two points.
The definition of $\ \Fl,\ \A,\ \L,\ \G\ $ is completely
analogous to the definition in Section 2.
The crucial step is to introduce an almost-grading and
to find dual systems of basis elements.
This is done in \cite\SCHLTH, \cite{\SCHLL, 3.ref.}.
For a quick review see \cite\SCHLCT. (In this context see also
Sadov \cite{\SADOV}
and the appendix of \cite\KNPMS.)
We want to recall here some steps.
Let $A$ be the finite set of points where poles are allowed
(which for $g\ge 1$ are in generic position).
The set $A$ has to be splited into  two non-empty disjoint
subset $I$ and $O$. The set $I$ is called the set of {\it in-points}, the
set $O$ the set of {\it out-points}.
Let $\# K$ be the number of in-points.
One fixes again a set of basis elements
$\ f_{n,p}^\l,\  n\in\Z,\; p=1,\ldots,K\ $
of $\Fl$ (the space of meromorphic forms of weight $\l$ which are
holomorphic on $\Gamma\setminus A$), by requiring certain
zero orders at the points of $A$.
To give  an example: Let
the number of in-points be equal to the number of out-points.
 $\l\ne 0,1$, and
$\  I=\{P_1,P_2,\ldots,P_K\}\ $,
$O=\{Q_1,Q_2,\ldots,Q_K\}\ $
be points in generic positions.
Then there is for every $n\in\Z$ and every $p=1,\ldots,K$ up to
multiplication with a scalar a  unique element
$f_{n,p}^\l\in\Fl$ with
$$
\aligned
\ord_{P_i}(f_{n,p}^\l)&=(n+1-\l)-\d_{i,p},\qquad i=1,\ldots,K,\\
\ord_{Q_i}(f_{n,p}^\l)&=-(n+1-\l),\qquad \qquad i=1,\ldots,K-1,\\
\ord_{Q_K}(f_{n,p}^\l)&=-(n+1-\l)+(2\l-1)(g-1)\ .
\endaligned\tag \kn-1
$$
This can be shown
either  by using Riemann-Roch type arguments
or by explicit constructions.
These elements will be the basis elements.
In the general situation there are modifications only at the out-points.
The element $f_{n,p}^\l$ is defined to be a homogeneous element of degree
$n$.
One fixes a differential $\rho\in\Fn 1$ which has exact pole order
1 at the points in $A$, positive residues at the points in $I$,
negative residues at the points in $O$, and purely
imaginary periods.  The level lines $C_\tau$ are defined completely
in the same manner as in (2-1).
Every level line separates the in- from
the out-points.
For $\tau\ll 0$ the level line $C_\tau$ is a disjoint union
of deformed circles around the points in $I$.
For $\tau\gg 0$ it is
a disjoint union
of deformed circles around the points in $O$.
\footnote{In the interpretation of string theory this
 $\tau$ might be interpreted
as proper time of the string on the world sheet.
}
We obtain the important duality
$$\cintt f_{n,p}^\l\cdot f_{m,r}^{1-\l}=\d_{n,-m}\cdot
\delta_{p,r}\ .\tag \kn-2$$
Again we use
$\ A_{n,p}=f_{n,p}^0,\  e_{n,p}=f_{n,p}^{-1},\
\omega^{n,p}=f^1_{-n,p},
\  \Omega^{n,p}=f^2_{-n,p}$.
One obtains an almost-graded structure  with respect to the
above introduced degree. Analogous formulas to (2-13) and (2-14)
are valid. For example
$$A_{n,p}\cdot A_{m,r}=\sum_{k=n+m}^{n+m+L}\sum_{s=1}^K
\a_{(n,p),(m,r)}^{(k,s)}
\,A_{k,s},\
\a_{(n,p),(m,r)}^{(k,s)}=\cintt A_{n,p}A_{m,r}\w^{k,s}\ .\tag\kn-3$$
Note that  we even have
$\ \a_{(n,p)(m,r)}^{(n+m,s)}=\delta_p^s\cdot\delta_r^s\ $.
Of course the constant $L$ (and $M$ in the equivalent formula to
(\knset-13) and all other introduced bounds) depend also
on the number of points in $A$ and their splitting into $I$ and $O$.
Explicit formulas can be found in \cite\SCHLTH.
The corresponding ``{\it delta distribution}'' is now
$$\Delta(Q',Q)=\sum_{n\in\Z}\sum_{p=1}^K A_{n,p}(Q')\w^{n,p}(Q)
\ .\tag\kn-4$$
For $\a\in H^1(\Gamma\setminus A,\Z)$
the 2-cocycles (2-3) and (2-6) define again central extensions
 $\Ah_\a$ and
$\Lh_\a$ of $\A$, resp\. of $\L$, and Equation (2-5) defines
a generalized multi-point affine Kac-Moody algebra  $\Gh_\a$.
If we choose a level line as integration cycle the cocycles will be local,
i.e\. we have similar bounds as in (2-15). Hence by defining
$\ \deg(t)=0\ $ we are able to extend the almost-grading to the centrally
extended algebras $\Ah,\ \Lh,\ \Gh\ $.
Again, without mentioning the cycle $\a$ we mean integration along a
level line $C_\tau$.

We have a decomposition
$\ \A=\A_-\oplus\A_0\oplus\A_+\ $, and $\
\G=\G_-\oplus\G_0\oplus\G_+
\ $
induced from $\  \A_-:=\langle A_n\mid n\le -P-1\rangle\ $,
$ \A_0:=\langle A_n\mid -P\le  n\le 0\rangle \ $,
$\ \A_+:=\langle A_n\mid n\ge 1\rangle\ $,
where $P$ is now a suitable positive constant.
Everything works as above. For $x\in\g$ we define $\ x(n,p):=x
\otimes A_{n,p}\in\Gh$. Hence the (non-central)
generators of $\Gh$ come with two labels.
Let $V$ be an admissible representation of $\Gh$ then we denote by
$x(A)$ the operator corresponding to the element $x\otimes A$. we use also
$x(n,p):=x(A_{n,p})$.

Where we had above a summation over $n$ we have to  add to this summation
another (finite) summation over $p=1,\ldots, K$:
$\ \widehat{u^i}(Q)=\sum_{n}\sum_p u^i(n,p)\cdot \w^{n,p}(Q)\ $.
Here and in the following the summation over the first index of the
label is always over $\Z$ and  over the second index of the label
over $1,\ldots, K$ if nothing else is said.
The higher genus (multi-point) {\it Sugawara  operator}
is defined as
$$T(Q):=\frac 12\sum_i:\widehat{u_i}(Q)\widehat{u^i}(Q):\quad =\
\frac 12\sum_{n,m}\sum_{p,s}\sum_i:u_i(n,p)u^i(m,s):
\w^{n,p}(Q)\w^{m,s}(Q)\ .\tag\kn-5$$
We decompose it again as
$$T(Q)=\sum_k\sum_r L_{k,r}\cdot\Omega^{k,r}(Q)\ ,\tag\kn-6$$
with
$$
\gathered
L_{k,r}=\cintt T(Q)e_{k,r}(Q)=\frac 12\sum_{n,m}\sum_{p,s}\sum_i
:u_i(n,p)u^i(m,s):l_{(k,r)}^{(n,p)(m,s)},\\
\text{where}\qquad
\ l_{(k,r)}^{(n,p)(m,s)}=\cintt w^{n,p}(Q)w^{m,s}(Q)e_{k,r}(Q)\ .
\endgathered\tag\kn-7$$
We define
$\
\gamma_{(n,p)(m,s)}:=\gamma(A_{n,p},A_{m,s}):=
\cintt  A_{n,p}(Q)dA_{m,s}(Q)\ $,
and obtain
\proclaim{Proposition \kn.1}
Let $\g$ be either an abelian or a simple Lie algebra, then
we have
$$\gather
[L_{k,r},x(n,p)]=-(\ce+\k)\sum_m\sum_s\,K_{(n,p),(k,r)}^{(m,s)} x(m,s)
\ ,\tag\kn-8
\\
[L_{k,r},x(n,p)]=-(\ce+\k)\,x(\nabla_{e_{k,r}}A_{n,p})\ ,\tag\kn-9
\\
[L_{k,r},\widehat{x}(Q)]=(\ce+\k)\,e_{k,r}\ldot \widehat{x}(Q)\ ,
\tag\kn-10
\\
\text{with}\quad
K_{(n,p),(k,r)}^{(m,s)}
=\cintt w^{m,s}e_{k,r}dA_{n,p}=
\sum_l\sum_v l_{(k,r)}^{(m,s)(l,v)}\gamma_{(l,v)(n,p)}\ .\tag\kn-11
\endgather
$$
The result does not depend on the normal ordering.
\endproclaim
If one checks  the proof in Section 4 one sees  that  essentially
only the duality and the  ``delta distribution'' has been used and
they have been generalized. Beside this one uses also
the generalization of the Lemma of Bonora et al.:
\proclaim{Lemma \kn.1}
For every $N$
$$
\bigg(\sum_{n>N\atop m\le N}-\sum_{m>N\atop n\le N}\bigg)
\sum_{p,s}A_{n,p}(Q')A_{m,s}(Q)\w^{m,s}(Q')\w^{n,p}(Q)
=d'\Delta(Q',Q)\ .\tag\kn-12
$$
Here $d'$ means differentiation with respect to the variable
$Q'$.
\endproclaim
\demo{Sketch of the  proof}
In the appendix of \cite\KNPMS\  Krichever and Novikov state the
uniqueness of a local 2-cocycle for the (multi-point) algebras
$\A$ (up to scalar) and $\L$ (up to scalar and coboundary).
Indeed the proof of the ``uniqueness'' presented in \cite\KNFA\ for
the 2-point situation using discrete Baker-Akhiezer functions
should  generalize  (see \cite{\SCHLTH, p.89}).
If we use this claim we can argue as in Section 4.
The right hand side of  (\knpr-11) in its suitable generalization has
to be
again a scalar multiple of  the cocycle (\knset-3) (with integration
over a level line).
 To fix the scalar  we calculate
$\ [A_{i,p},A_{-i,p}]\,\Phi$. Note that
we have  $A_{k,s}A_{n,r}=A_{k+n,r}\delta_s^r+$
$A$-terms with higher degree as $k+n$ and
$\gamma(A_{i,p},A_{-i,p})=(-i)$
(there is only a residue at the point $P_p$).
If we want to avoid the statement about the uniqueness we
have to make lengthy local calculations as indicated in \cite\BRRW.
\qed
\enddemo
Now Proposition 3.3
and its proof have their obvious generalizations
\proclaim{Proposition \kn.2}
The operators $L_k\in gl(V)$ and $id=1\in gl(V)$ close up to
a Lie subalgebra of $gl(V)$ with the commutator relation
$$[L_{k,r},L_{l,v}]=-(\ce+\k)\sum_n\sum_p C_{(k,r)(l,v)}^{(n,p)}L_{n,p}
-\frac 12
\ce(\ce+\k)\dim \g\cdot \chi_{(k,r)(l,v)}\cdot id
\ .\tag\kn-13$$
$\chi_{(k,r)(l,v)}\ne 0$ implies that
$  -P\le k+l\le 0$ with a positive constant $P$ not depending on
 $k$ and $l$.
\endproclaim
\noindent
The splitting of $\chi_{(k,r)(l,s)}$ and its expression is
completely analogous to (\knsu-11).
The same is true for the expression (\knsu-12) for the form of
the cocycle for certain basis elements.

We get the
\proclaim{Theorem \kn.1}
Let $\g$ be either a finite dimensional abelian or simple Lie algebra and
 $2\k$ be the eigenvalue of the Casimir operator in the adjoint
representation.
Let $A$ be a finite set of points (for $g\ge 1$ in general position)
and $A=I\cup O$ a splitting into two disjoint non-empty subsets.
Let $\A$ be the algebra of meromorphic functions which are
holomorphic on $\Gamma\setminus A$ and $\Gh$
be the generalized higher genus multi-point affine Kac-Moody algebra, i.e\.
the
central extension of $\G=\g\otimes \A$ defined  by
$$[\widehat{x\otimes f},\widehat{y\otimes g}]=
\widehat{[x,y]\otimes (f g)}-(x|y)\cdot\cintt fdg\cdot t,\qquad
[\,t,\Gh_\a]=0\ ,\tag\kn-14$$
where $C_\tau$ is a level line separating the points in $I$ and $O$.
Let $V$ be  an admissible representation  where the central element
$t$ operates as  $\ce\cdot identity$. If $\ce+\k\ne 0$ then the rescaled
modes
$$L_{k,r}^*=\frac {-1}{2(\ce+\k)}\sum_{n,m}\sum_{p,s}\sum_i
 :u_i(n,p)u^i(m,s):l_{(k,r)}^{(n,p)(m,s)}
\ ,$$
of the Sugawara operator define a representation of
a local central extension $\Lh'$ of the Krichever Novikov vector field
algebra $\L$.
\endproclaim
Using again the statement of Krichever and  Novikov
about the uniqueness of a local cocycle,
and the corresponding result to Lemma~\knsu.2 (see\cite{\SCHLTH})
we get indeed
a representation of $\Lh$ defined by the cocycle (\knsu-15).
%
%
%
{
}
\enddocument
\bye